\pdfoutput=1
\RequirePackage{ifpdf}
\ifpdf 
\documentclass[pdftex]{sigma}
\else
\documentclass{sigma}
\fi

\numberwithin{equation}{section}

\newtheorem{Theorem}{Theorem}[section]
\newtheorem{Corollary}[Theorem]{Corollary}
\newtheorem{Lemma}[Theorem]{Lemma}
\newtheorem{Proposition}[Theorem]{Proposition}
 { \theoremstyle{definition}
\newtheorem{Definition}[Theorem]{Definition}
\newtheorem{Remark}[Theorem]{Remark}
}

\usepackage{mathrsfs}
\usepackage{enumitem}

\newcommand{\standout}[1]{\textit{#1}} 

\newcommand{\n}[1]{{#1}_{[n]}}

\newcommand{\actby}[1]{{}^{#1}}
\newcommand{\actbylessspace}[1]{{}^{#1}\!}

\newcommand{\vol}[1]{\text{vol}_{#1}^{\perp}}
\newcommand{\HHSSG}[1]{\HH^{#1}(S(V),S(V)\#G)}

\newcommand{\kappaLxy}[5]{\kappa^{*}_{#1}\big(#3+{}^{#2}#3,\kappa^L_{#2}(#4,#5)\big)}
\newcommand{\phibulletxy}[5]{\kappaLxy{#1}{#2}{#3}{#4}{#5}+\kappaLxy{#1}{#2}{#4}{#5}{#3}+\kappaLxy{#1}{#2}{#5}{#3}{#4}}

\DeclareMathOperator{\HH}{HH}
\DeclareMathOperator{\Ext}{Ext}
\DeclareMathOperator{\Hom}{Hom}
\DeclareMathOperator{\Span}{Span}
\DeclareMathOperator{\Sym}{Sym}
\DeclareMathOperator{\Alt}{Alt}
\DeclareMathOperator{\im}{im}
\DeclareMathOperator{\codim}{codim}
\DeclareMathOperator{\gr}{gr}
\DeclareMathOperator{\refl}{ref}
\DeclareMathOperator{\tri}{3-cyc}

\begin{document}

\newcommand{\arXivNumber}{1912.06743}

\renewcommand{\PaperNumber}{039}

\FirstPageHeading

\ShortArticleName{Degree-One Rational Cherednik Algebras for the Symmetric Group}

\ArticleName{Degree-One Rational Cherednik Algebras\\ for the Symmetric Group}

\Author{Briana FOSTER-GREENWOOD~$^{\rm a}$ and Cathy KRILOFF~$^{\rm b}$}

\AuthorNameForHeading{B.~Foster-Greenwood and C.~Kriloff}

\Address{$^{\rm a)}$~Department of Mathematics and Statistics, California State Polytechnic University,\\
\hphantom{$^{\rm a)}$}~Pomona, California 91768, USA}
\EmailD{\href{mailto:brianaf@cpp.edu}{brianaf@cpp.edu}}

\Address{$^{\rm b)}$~Department of Mathematics and Statistics, Idaho State University,\\
\hphantom{$^{\rm b)}$}~Pocatello, Idaho 83209, USA}
\EmailD{\href{mailto:cathykriloff@isu.edu}{cathykriloff@isu.edu}}

\ArticleDates{Received August 07, 2020, in final form April 02, 2021; Published online April 19, 2021}

\Abstract{Drinfeld orbifold algebras deform skew group algebras in polynomial degree at~most one and hence encompass graded Hecke algebras, and in particular symplectic reflection algebras and rational Cherednik algebras. We~introduce parametrized families of~Drinfeld orbifold algebras for symmetric groups acting on doubled representations that generalize rational Cherednik algebras by deforming in degree one. We~characterize rich fami\-lies of~maps recording commutator relations with their linear parts supported only on and only off the identity when the symmetric group acts on the natural permutation representation plus its dual. This produces degree-one versions of $\mathfrak{gl}_n$-type rational Cherednik algebras. When the symmetric group acts on the standard irreducible reflection representation plus its dual there are no degree-one Lie orbifold algebra maps, but there is a~three-parameter family of~Drinfeld orbifold algebras arising from maps supported only off the identity. These provide degree-one generalizations of the $\mathfrak{sl}_n$-type rational Cherednik algebras $H_{0,c}$.}


\Classification{16S80; 16E40; 16S35; 20B30}

\section{Introduction}

Skew group or smash-product algebras $S(V)\# G$ twist the symmetric algebra $S(V)$ of a
finite-dimensional vector space $V$ together with the action of the group algebra $\mathbb{C} G$
of a finite group $G$ acting linearly on $V$. The center is the invariant polynomial ring
$S(V)^G$ and there is a natural grading by polynomial degree, with elements in~$V$ of
degree one and elements in $\mathbb{C} G$ of degree zero.

Utilizing parameter maps that originate as Hochschild $2$-cocycles to explore formal
deformations of $S(V) \# G$ has proven useful because although the resulting algebras
are noncommutative they give rise to deformations of $S(V)^G$ (by examining centers),
yet are easily described
as quotient algebras. Both the polynomial degree and the support, i.e., which group elements
appear in the nonzero image, are helpful descriptors for the parameter maps and hence
the relations for~the quotient algebras.

Degree-zero deformations of skew group algebras involve parameter maps that
identify commutators of elements in~$V$ with certain elements of the
group algebra. Several important families of these are of broad
interest in noncommutative geometry, combinatorics, and representation theory
and are the subject of an already extensive literature (see~\cite{Gordon2008}
and~\cite{Gordon2010} and further references therein). By~comparison, finding elements of degree one with which to identify commutators
of elements in~$V$ requires a more intricate analysis of which cocycles pass obstructions in
cohomology in order to determine if the resulting deformations satisfy PBW properties~\cite{FGK, SWorbifold}.
As a result, degree-one deformations are not as well understood or as often studied, yet
could also be significant in giving insight into deformations of the invariant algebra
$S(V)^G$ and in~con\-nection with singularities of orbifolds.

Degree-zero deformations of skew group algebras are called Drinfeld graded Hecke algebras in~recognition of their origins in~\cite{Drinfeld1986} (see also~\cite{Lusztig1988}).
These include the important special cases when $G$ acts on a symplectic vector space,
and more particularly when $G$ is a complex reflection group acting by the sum of a
reflection representation and its dual (a doubled representation). The~lat\-ter leads
to the rational Cherednik algebras, first introduced in~\cite{Cherednik1991} as rational
dege\-ne\-ra\-tions of~double affine Hecke algebras and later highlighted as an important
subfamily of~the more general symplectic reflection algebras introduced in~\cite{EtingofGinzburg2002}.
When built from an action of~the symmetric group, rational Cherednik
algebras model Hamiltonian reduction in quantum mechanics and are used to show
integrability of Calogero--Moser systems~\cite{Etingof2007}.

Degree-one deformations of skew group algebras were termed Drinfeld orbifold algebras and
characterized via explicit PBW conditions on parameter maps in~\cite{SWorbifold}, building on~\cite{Bergman1978} and~\cite{BravermanGaitsgory}.
The conditions are also interpreted in Hochschild cohomology.
In~\cite{FGK} we describe the Drinfeld orbifold algebras for $S_n$ acting on its natural
permutation representation, $W\cong \mathbb{C}^n$, by starting with candidate $2$-cocycles and imposing the
PBW conditions from~\cite{SWorbifold}.
Here we expand on that class of examples by considering $S_n$
acting on both its doubled permutation representation $W^*\oplus W$ and the doubled representation
$\mathfrak{h}^*\oplus \mathfrak{h}$, where $W=\mathfrak{h} \oplus \iota$ is the sum of the $(n-1)$-dimensional irreducible
standard and the trivial representations. This not only results in much richer families
of algebras, but also yields degree-one generalizations of rational Cherednik algebras for these
doubled representations.

More specifically, in~\cite{FGK} we describe all Drinfeld orbifold algebras
where the linear parts of~the maps recording commutator relations are supported only
on or only off the identity in~$S_n$, and show there are no such maps with
linear part supported both on and off the identity. For~the two
doubled representations of $S_n$ considered here we describe all degree-one
families of~Drinfeld orbifold algebras whose maps have linear part supported only on
or only off the identity (Theorems~\ref{Lieorbifold} and~\ref{permexamples}). For~maps with linear part supported both on and off the identity
we provide a family of examples involving $W^*\oplus W$
(Theorem~\ref{thm:DOAMapsCombined})
and observe there are no corresponding such maps for
the doubled standard representation $\mathfrak{h}^*\oplus \mathfrak{h}$
(Remark~\ref{NoExtensionsOfExamples}). We~summarize our main results.
\begin{Theorem}\label{thA}
For the symmetric group $S_n$ ($n\geq 3)$ acting on $V \cong \mathbb{C}^{2n}$
by the doubled permutation representation, there is
\begin{enumerate}\itemsep=0pt
\item[$(1)$] a $17$-parameter family of Lie orbifold algebras described by $22$
homogeneous quadratic equations, and
\item[$(2)$] a seven-parameter family of~Drinfeld orbifold algebras described in terms
of parameter maps with linear part supported only off the identity that are
controlled by four homogeneous quadratic equations in six of the parameters.
\end{enumerate}
These are the only degree-one deformations of the skew group algebra for $S_n$
acting by the doubled permutation representation whose
parameter maps have linear part supported only on or only off the identity.
\end{Theorem}
See Theorems~\ref{thm:LOAMapsDoubledPerm} and~\ref{thm:DOAMapsDoubledPerm}
for more details about the maps, Theorems~\ref{Lieorbifold} and~\ref{permexamples}
for the resulting quotient algebras, and Table~\ref{SummaryDOAMapsDoubledPerm}
in~Section~\ref{AlgsDoubledPermRep} for a summary.

\begin{Theorem}\label{thB}For the symmetric group $S_n$ $(n\geq 3)$ acting on $V \cong \mathbb{C}^{2n-2}$ by the
doubled standard representation, there are no degree-one Lie orbifold algebras,
but there is a three-parameter family of~Drinfeld orbifold algebras described by
parameter maps with linear part supported only off the identity.

These are the only degree-one deformations of the skew group algebra for $S_n$
acting by the doubled standard representation whose
parameter maps have linear part supported only on or only off the identity.
\end{Theorem}
See Theorems~\ref{thm:LOAMapsDoubledStd} and~\ref{thm:DOAMapsDoubledStd}
for details, Theorem~\ref{refexamples} for the resulting algebras,
and Table~\ref{SummaryDOAMapsDoubledStd} in~Section~\ref{AlgsDoubledStdRep} for a~summary.
The $S_2$ case in Theorems~\ref{thA} and~\ref{thB} can be analyzed in a similar way but there are some
differences in the dimensions of spaces of pre-Drinfeld orbifold algebra maps and
in the explicit PBW conditions.

The algebras in Theorems~\ref{thA} and~\ref{thB} specialize
to the well-known rational Cherednik algebras for the symmetric group,
described as of $\mathfrak{gl}_n$- and $\mathfrak{sl}_n$-type respectively
in~\cite{GanGinzburg2006}, and hence should be of substantial interest.
In particular, Theorem~\ref{thB}
provides a degree one version of~the~$\mathfrak{sl}_n$-type rational Cherednik
algebras $H_{0,c}$. We~refer to the algebras as degree-one rational Cherednik algebras.
Investigating their structure, properties, combinatorics, representation theory,
geometric significance, and potential importance in physics
should provide fertile ground for future research. It would also be natural
to explore whether similar algebras exist for other complex reflection groups.

The paper is organized as follows.
After a brief summary of preliminaries in~Section~\ref{Preliminaries}
that apply to any finite group acting linearly on $\mathbb{C}^n$, we restrict
to the setting of the symmetric group and the two doubled representations
of interest, except as noted in Lemmas~\ref{eigenvectors} and~\ref{le:orbitphixyphig}, Proposition~\ref{prop:directsum},
and Corollary~\ref{prop:Extension}.
All pre-Drinfeld orbifold algebra maps for $S_n$ acting by the doubled
permutation representation are constructed in~Section~\ref{InvSkewSymmBilMaps}. We~analyze when these lift in~Sections~\ref{InvLieBrackets} and~\ref{Lifting}, proving Theorems~\ref{thm:LOAMapsDoubledPerm},~\ref{thm:DOAMapsDoubledPerm},
and~\ref{thm:DOAMapsCombined} using computational details treated
earlier in the two sections.
In particular, Sections~\ref{AllVecInW}--\ref{TwoVecInWTwoIndices} provide
explicit equations governing the parameter maps described in
Theorem~\ref{thm:LOAMapsDoubledPerm} and Section~\ref{algvars} provides some related
algebraic varieties that may be of independent interest.
Section~\ref{DoubledStandardDefs} begins with Proposition~\ref{prop:directsum} providing conditions under which we can combine Drinfeld orbifold algebra maps for subrepresentations into a map for their direct sum.
Corollary~\ref{prop:Extension} is then used with the results from
Sections~\ref{InvLieBrackets} and~\ref{Lifting} to describe in Theorems~\ref{thm:LOAMapsDoubledStd}
and~\ref{thm:DOAMapsDoubledStd} all Drinfeld orbifold algebra maps
for $S_n$ acting by the doubled standard representation on the subspace
$\mathfrak{h}^* \oplus \mathfrak{h}$ when the linear part is supported only on or only off
the identity. In Section~\ref{Deg1RCAs} we present as quotients
the resulting degree-one rational Cherednik algebras arising from the maps in~Sections~\ref{InvLieBrackets}--\ref{DoubledStandardDefs}.

\section{Preliminaries}\label{Preliminaries}
Throughout this section, we let $G$ be a finite group acting linearly on a vector space $V\cong\mathbb{C}^n$. All tensors will be over $\mathbb{C}$.

\subsection{Skew group algebras}
Let $G$ be a finite group that acts on a $\mathbb{C}$-algebra $R$ by
algebra automorphisms, and write $\actby{g}s$ for the result of acting by $g\in G$ on $s\in R$. The \standout{skew group algebra} $R\#G$ is the
semi-direct product algebra $R\rtimes \mathbb{C} G$ with underlying vector space $R\otimes \mathbb{C} G$ and multiplication of simple tensors defined by
\begin{gather*}
(r\otimes g)(s\otimes h)=r(\actby{g}s)\otimes gh
\end{gather*}
for all $r,s\in R$ and $g,h\in G$.
The skew group algebra becomes a $G$-module by letting $G$ act diagonally on
$R\otimes \mathbb{C} G$, with conjugation on the group algebra factor:
\begin{gather*}
\actby{g}(s\otimes h)=(\actby{g}s)\otimes(\actby{g}h)=(\actby{g}s)\otimes ghg^{-1}.
\end{gather*}
In working with elements of skew group algebras, we commonly omit tensor symbols unless
the tensor factors are lengthy expressions.

If $G$ acts linearly on
a vector space $V\cong\mathbb{C}^n$, then $G$ also acts on the tensor algebra $T(V)$ and
symmetric algebra $S(V)$ by algebra automorphisms.
Assign elements of $V$ degree one and elements of $G$ degree zero
to make the skew group algebras $T(V)\# G$ and $S(V)\# G$ graded algebras.

\subsection{Cochains}

A \standout{$k$-cochain} is a $G$-graded linear map $\mu=\sum_{g\in G}\mu_g g$ with components $\mu_g\colon \bigwedge^k V\to S(V)$. If~each $\mu_g$ maps into $V$, then $\mu$ is called a
\standout{linear cochain}, and if each $\mu_g$ maps into $\mathbb{C}$, then $\mu$ is called
a \standout{constant cochain}.

We regard a map $\mu$ on $\bigwedge^kV$ as a multilinear alternating map on $V^k$ and write $\mu(v_1,\ldots,v_k)$ in place of $\mu(v_1\wedge\cdots\wedge v_k)$.
Of course, if
$\mu(v_1,\ldots,v_k)=0$, then $\mu$ is zero on any permutation of
$v_1,\ldots,v_k$. Also, if $\mu$ is zero on all $k$-tuples of basis vectors, then $\mu$ is zero on any $k$-tuple of~vectors. We~exploit these facts in the computations in~Sections~\ref{InvLieBrackets} and~\ref{Lifting}.

The \standout{support of a cochain $\mu$} is the set of group elements for which the component $\mu_g$ is not the zero map. For~$X$ a subset of $G$, we say a cochain $\mu$ is \standout{supported only on {\boldmath $X$}} if $\mu_g=0$ for all~$g$ not in $X$. Similarly, we say $\mu$ is
\standout{supported only off {\boldmath $X$}} if $\mu_g=0$ for all $g$ in $X$.
At times, it is convenient to talk about support in a weaker sense, so we say $\mu$ is \standout{supported on {\boldmath $X$}} if $\mu_g\neq0$ for some $g$ in $X$
and that $\mu$ is \standout{supported off {\boldmath $X$}} if $\mu_g\neq0$ for some $g$ not in $X$. (Hence, it is possible for a cochain to be simultaneously supported on and off of a set.)
The \standout{kernel of a cochain $\mu$} is the set of vectors $v_0$ such that $\mu(v_0,v_1,\ldots,v_{k-1})=0$ for all $v_1,\ldots,v_{k-1}\in V$.

The group $G$ acts on the components of a cochain. Specifically, for a group element $h$ and component $\mu_g$, the map
$\actby{h}\mu_g$ is defined by $\big(\actby{h}\mu_g\big)(v_1,\ldots,v_k)=\actby{h}\big(\mu_g\big(\actby{h^{-1}}v_1,\ldots,\actby{h^{-1}}v_k\big)\big)$. In turn, the group acts on the space of cochains by letting $\actby{h}\mu=\sum_{g\in G}\actby{h}\mu_g\otimes hgh^{-1}$. Thus $\mu$ is a~\standout{$G$-invariant cochain} if and only if $\actby{h}\mu_g=\mu_{hgh^{-1}}$ for all $g,h\in G$.

\subsection{Drinfeld orbifold algebras}
For a parameter map $\kappa=\kappa^L+\kappa^C$, where $\kappa^L$ is a linear $2$-cochain and $\kappa^C$ is a constant $2$-cochain, the quotient algebra
\begin{gather*}
\mathcal{H}_{\kappa}=T(V)\#G/\left\langle vw-wv-\kappa^L(v,w)-\kappa^C(v,w)\mid v,w\in V\right\rangle
\end{gather*}
is called a \standout{Drinfeld orbifold algebra} if the associated graded algebra $\gr\mathcal{H}_{\kappa}$ is isomorphic to the skew group algebra
$S(V)\# G$. The condition $\gr\mathcal{H}_{\kappa}\cong S(V)\# G$ is
called a \standout{Poincar\'e--Birkhoff--Witt $($PBW$)$ condition}, in analogy with
the PBW Theorem for universal enveloping algebras.

Further, if $\mathcal{H}_{\kappa}$ is a Drinfeld orbifold algebra and $t$ is a complex parameter,
then
\begin{gather*}
\mathcal{H}_{\kappa,t}:=T(V)\#G[t]/\left\langle vw-wv-\kappa^L(v,w)t-\kappa^C(v,w)t^2\mid v,w\in V\right\rangle
\end{gather*}
is called a \standout{Drinfeld orbifold algebra over $\mathbb{C}[t]$}.
In~\cite[Theorem~2.1]{SWorbifold}, Shepler and Witherspoon make an explicit
connection between the PBW condition and deformations in the sense of~Ger\-s\-tenhaber~\cite{GerstenhaberSchack} by showing how to interpret Drinfeld orbifold
algebras over $\mathbb{C}[t]$ as formal defor\-ma\-tions of the skew group algebra
$S(V)\# G$. For~more on the broader context of formal deformations
see~\cite[Section~4]{FGK}.

\subsection{Lie orbifold algebras}
The parameter maps of~Drinfeld orbifold algebras decompose as
$\kappa=\sum_{g} \kappa_g g$. When $\kappa$ is a~para\-meter map for a Drinfeld orbifold algebra and the linear part $\kappa^L=\kappa^L_1$
is supported only on the identity then the map gives rise to a Lie orbifold algebra
(see~\cite[Section~4]{SWorbifold} and~Definition~\ref{def:fourconditions}).
Lie orbifold algebras deform universal enveloping
algebras twisted by a group action just as certain symplectic reflection algebras
deform Weyl algebras twisted by a group action.

\subsection{Drinfeld orbifold algebra maps}
Though the defining PBW condition for a Drinfeld orbifold algebra $\mathcal{H}_{\kappa}$ involves an isomorphism of algebras, Shepler and Witherspoon proved an equivalent characterization~\cite[Theorem~3.1]{SWorbifold} in terms of properties of the parameter map $\kappa$.

\begin{Definition}\label{def:fourconditions}
	Let $\kappa=\kappa^L+\kappa^C$ where $\kappa^L$ is a linear $2$-cochain and $\kappa^C$ is a constant $2$-cochain, and let $\Alt_3$ denote the alternating group on three elements. Let $V^g$ denote the set of vectors in~$V$ that are fixed by group element $g$.
	We say $\kappa$ is a {\it Drinfeld orbifold algebra map} if the following conditions are satisfied for all $g\in G$ and $v_1,v_2,v_3\in V$:
\begin{gather}
\label{PBWconditions-0}
\im\kappa^L_g\subseteq V^g,
\\
\label{PBWconditions-i}
\text{the map $\kappa$ is $G$-invariant,}
\\
\label{PBWconditions-ii} \sum_{\sigma\in\Alt_3}\kappa^L_g(v_{\sigma(2)},v_{\sigma(3)})(\actby{g}v_{\sigma(1)}-v_{\sigma(1)})=0\qquad
\text{in}\quad S(V),
\\
\label{PBWconditions-iii}
\sum_{\sigma\in\Alt_3}\sum_{xy=g}
\kappa^L_x\big(v_{\sigma(1)}+\actby{y}v_{\sigma(1)},\kappa^L_y(v_{\sigma(2)},v_{\sigma(3)})\big) =2\!\!\sum_{\sigma\in\Alt_3}\!\kappa^C_g(v_{\sigma(2)},v_{\sigma(3)})(\actby{g}v_{\sigma(1)}-v_{\sigma(1)}),
\\
\label{PBWconditions-iv}
\sum_{\sigma\in\Alt_3}\sum_{xy=g}
\kappa^C_x\big(v_{\sigma(1)}+\actby{y}v_{\sigma(1)},\kappa^L_y(v_{\sigma(2)},v_{\sigma(3)})\big)=0.
\end{gather}
In the special case when the linear component $\kappa^L$ of	a Drinfeld orbifold
algebra map is supported only on the identity, we call $\kappa$ a {\it Lie orbifold algebra map}.
\end{Definition}

To simplify reference to the expressions appearing in the last three Drinfeld orbifold algebra map properties, we define operators $\psi$ and $\phi$ that convert $2$-cochains (such as $\kappa^L$ and $\kappa^C$) into the $3$-cochains we see evaluated within the properties.

\begin{Definition}
	Let $\mu$ denote a linear or constant $2$-cochain and $\nu$ a linear $2$-cochain.
	Define $\psi(\mu)=\sum_{g\in G}\psi_gg$ to be the $3$-cochain with components $\psi_g\colon \bigwedge^3V\rightarrow S(V)$ given by
	\begin{gather*}
\psi_g(v_1,v_2,v_3)=\sum_{\sigma\in\Alt_3}\mu_g(v_{\sigma(1)},v_{\sigma(2)})(\actby{g}v_{\sigma(3)}-v_{\sigma(3)}).
	\end{gather*}
	Define $\phi(\mu,\nu)=\sum_{g\in G}\phi_g g$ to be the $3$-cochain
	with components $\phi_g=\sum_{xy=g}\phi_{x,y}$, where $\phi_{x,y}\colon \bigwedge^3 V\to V\oplus\mathbb{C}$ is given by
	\begin{equation}
	\label{phi}
	\phi_{x,y}(v_1,v_2,v_3)=\sum_{\sigma\in\Alt_3}\mu_x(v_{\sigma(1)}+\actby{y}v_{\sigma(1)},\nu_y(v_{\sigma(2)},v_{\sigma(3)})).
	\end{equation}
	Thus $\phi(\mu,\nu)$ is $G$-graded with components $\phi_g$ and also $(G\times G)$-graded with components $\phi_{x,y}$.
\end{Definition}

For the interested reader, we indicate in~\cite{FGK} how the maps $\psi$ and $\phi$ relate to coboundary and bracket operations in Hochschild cohomology of a skew group algebra.

\subsection{Drinfeld orbifold algebra maps (condensed definition)}

Equipped with the definitions of $\psi$ and $\phi$, the properties of a Drinfeld orbifold map $\kappa=\kappa^L+\kappa^C$
(Definition~\ref{def:fourconditions}) may be expressed succinctly:
\begin{enumerate}	\itemsep=0pt
	\item[\eqref{PBWconditions-0}] \ $\im\kappa^L_g\subseteq V^g$ for each $g$ in $G$,
	\item[\eqref{PBWconditions-i}] \ the map $\kappa$ is $G$-invariant,
	\item[\eqref{PBWconditions-ii}] \ $\psi\big(\kappa^L\big)=0$,
	\item[\eqref{PBWconditions-iii}] \ $\phi\big(\kappa^L,\kappa^L\big)=2\psi\big(\kappa^C\big)$,
	\item[\eqref{PBWconditions-iv}] \ $\phi\big(\kappa^C,\kappa^L\big)=0$.
\end{enumerate}

Note that any $G$-invariant 2-cochain whose linear part is supported only on the identity
trivially satisfies properties~\eqref{PBWconditions-0}
and~\eqref{PBWconditions-ii}, so in this case it is enough to analyze
conditions under which properties~\eqref{PBWconditions-iii} and~\eqref{PBWconditions-iv} hold
(see Theorem~\ref{thm:LOAMapsDoubledPerm} and Section~\ref{InvLieBrackets}).

\begin{Remark}
	If $\mathcal{H}_{\kappa}$ is a Drinfeld orbifold algebra, then $\kappa$ must satisfy conditions~\eqref{PBWconditions-i}--\eqref{PBWconditions-iv}, but not necessarily the image constraint~\eqref{PBWconditions-0}.
	However,~\cite[Theorem~7.2(ii)]{SWorbifold} guarantees there will exist a
	Drinfeld orbifold algebra $\mathcal{H}_{\widetilde{\kappa}}$ such that
	$\mathcal{H}_{\widetilde{\kappa}}\cong\mathcal{H}_{\kappa}$ as
	filtered algebras and $\widetilde{\kappa}$ satisfies the image constraint
	$\im \widetilde{\kappa}^L_g\subseteq V^g$ for each $g$ in~$G$.
	Thus, in classifying Drinfeld orbifold algebras, it suffices to only
	consider Drinfeld orbifold algebra maps.
\end{Remark}

\begin{Theorem}[{\cite[Theorems~3.1 and~7.2(ii)]{SWorbifold}}]
	A quotient algebra $\mathcal{H}_{\kappa}$ satisfies the PBW condition
	$\gr\mathcal{H}_{\kappa}\cong S(V)\# G$ if and only if there exists a
	Drinfeld orbifold algebra map $\widetilde{\kappa}$ such that $\mathcal{H_{\kappa}}\cong\mathcal{H_{\widetilde{\kappa}}}$.
\end{Theorem}

\subsection{Strategy}
As described and utilized in~\cite{FGK}, the process of determining the set
of all Drinfeld orbifold algebra maps consists of two phases, and language
from cohomology and deformation theory can be used to describe each phase.
First, one finds all \standout{pre-Drinfeld orbifold algebra maps},
i.e., all $G$-invariant linear $2$-cochains satisfying the image
condition~\eqref{PBWconditions-0} and the mixed Jacobi identity~\eqref{PBWconditions-ii}.
To~find such maps supported on transpositions (Proposition~\ref{prop:kapparefl}) we utilize
a~bijec\-tion between pre-Drinfeld orbifold algebra maps and a particular
set of representatives of~Hochschild cohomology classes (see Lemma~\ref{ineedaname}).
But to find such maps supported only on~the identity (Proposition~\ref{prop:kappa1})
we present a simpler argument based on Lemma~\ref{eigenvectors}
analyzing the eigenvector structure of~images dependent on the group action
on input vectors.
In the second phase (Sections~\ref{InvLieBrackets} and~\ref{Lifting})
we determine for which
pre-Drinfeld orbifold algebra maps $\kappa^L$ there exists a compatible $G$-invariant
constant $2$-cochain $\kappa^C$ such that properties~\eqref{PBWconditions-iii}
and \eqref{PBWconditions-iv} hold. We~say $\kappa^C$
\standout{clears the first obstruction} if property~\eqref{PBWconditions-iii} holds
and \standout{clears the second obstruction} if property~\eqref{PBWconditions-iv} holds. If~a $G$-invariant constant $2$-cochain $\kappa^C$ clears both obstructions, then
we say $\kappa^L$ \standout{lifts} to the Drinfeld orbifold algebra map $\kappa=\kappa^L+\kappa^C$.

\subsection{Hochschild cohomology to pre-DOA maps}

We briefly recall how Hochschild cohomology can be used in general
to find linear and constant $2$-cochains $\kappa$ that are both $G$-invariant
and satisfy property~\eqref{PBWconditions-ii}. For~more detailed background discussion about the connections to deformations
and further references see~\cite{FGK}.

For an algebra $A$ over $\mathbb{C}$ with bimodule $M$,
the \standout{Hochschild cohomology} of $A$ with coefficients in $M$ is
$\HH^{\bullet}(A,M):=\Ext_{A\otimes A^{\text{op}}}^{\bullet}(A,M)$, which is
abbreviated to $\HH^{\bullet}(A)$ if $M=A$. For~any finite group $G$ acting linearly on a
vector space $V \cong \mathbb{C}^n$, and for $A=S(V)\# G$,
using results of~\c{S}te\-fan~\cite{Stefan} yields the following,
where $R^G$ denotes the set of elements in $R$ fixed by every $g$ in~$G$,
\begin{gather*}
\HH^{\bullet}(S(V)\#G) \cong\HHSSG{\bullet}^G \cong (H^{\bullet})^G.
\end{gather*}
Here $H^{\bullet}$ is the $G$-graded vector space
$H^{\bullet}=\bigoplus_{g\in G}H^{\bullet}_g$
with components
\begin{gather*}
H^{p,d}_g=S^d(V^{g}) \otimes\bigwedge^{p-\codim(V^g)}(V^g)^{*} \otimes
\bigwedge^{\codim(V^g)}\bigl((V^g)^*\bigr)^{\perp} \otimes
\mathbb{C}g,
\end{gather*}
first described independently by Farinati~\cite{Farinati} and
Ginzburg--Kaledin~\cite{GinzburgKaledin}.
Note that $H^{\bullet}$ is tri-graded by cohomological degree $p$,
homogeneous polynomial degree $d$, and group element $g$.

Since the exterior factors of $H_g^{p,d}$ can be identified
with a subspace of $\bigwedge^pV^*$, and since
$S^d(V^g)\otimes\bigwedge^pV^*\otimes\mathbb{C} g\cong\Hom\big(\bigwedge^p V,S^d(V^g)g\big)$,
the space $H^{\bullet}$ may be identified with a subspace of~the
cochains introduced earlier in this section. The next lemma
records the relationship bet\-ween properties~\eqref{PBWconditions-i}
and~\eqref{PBWconditions-ii} of a Drinfeld orbifold algebra map and Hochschild cohomology.
When $d=1$, the lemma is a restatement of~\cite[Theorem~7.2(i) and (ii)]{SWorbifold}.
When $d=0$, the lemma is a restatement of~\cite[Corollary 8.17(ii)]{SheplerWitherspoon08}.
It is also possible to give a linear algebraic proof 
in the spirit of~\cite[Lemma~1.8]{RamShepler}.

\begin{Lemma}\label{ineedaname}
	For a $2$-cochain $\kappa=\sum_{g\in G}\kappa_g g$ with $\im \kappa_g\subseteq S^d(V^g)$ for each $g\in G$,
	the following are equivalent:
	\begin{enumerate}\itemsep=0pt
		\item[$(a)$] The map $\kappa$ is $G$-invariant and satisfies
		the mixed Jacobi identity, i.e., for all $v_1,v_2,v_3\in V$
		\begin{gather*}
		[v_1,\kappa(v_2,v_3)]+[v_2,\kappa(v_3,v_1)]+[v_3,\kappa(v_1,v_2)]=0\qquad
\text{in}\quad S(V)\# G,
		\end{gather*}
		where $[\cdot,\cdot]$ denotes the commutator in $S(V)\# G$.
		\item[$(b)$] For all $g,h\in G$ and $v_1,v_2,v_3\in V$:
		\begin{enumerate}\itemsep=0pt
			\item[$(i)$] $\actby{h}(\kappa_g(v_1,v_2))=\kappa_{hgh^{-1}}\big(\actby{h}v_1,\actby{h}v_2\big)$, and
			\item[$(ii)$] $\kappa_g(v_1,v_2)(\actby{g}v_3-v_3)
			+\kappa_g(v_2,v_3)(\actby{g}v_1-v_1)
			+\kappa_g(v_3,v_1)(\actby{g}v_2-v_2)=0$.
		\end{enumerate}
		\item[$(c)$] The map $\kappa$ is an element of
		\begin{gather*}
		(H^{2,d})^G=\bigg(\bigoplus_{g \in G}
		\Bigl(
		S^d(V^{g})g \otimes
		\bigwedge^{2-\codim(V^g)}(V^g)^{*} \otimes
		\bigwedge^{\codim(V^g)}\bigl((V^g)^*\bigr)^{\perp}
		\Bigr)\bigg)^G.
		\end{gather*}
	\end{enumerate}
\end{Lemma}

\begin{Remark}\label{remark:codim}
	Part ($b(ii)$) of Lemma~\ref{ineedaname} is $2\psi(\kappa)=0$.
	Part ($c$) of Lemma~\ref{ineedaname} shows that $\kappa$ can only be supported
	on elements $g$ with $\codim V^g\in\{0,2\}$ since negative exterior powers
	are zero and an element $g$ with codimension one acts nontrivially on $H_g^{2,d}$.
\end{Remark}

\section{Pre-Drinfeld orbifold algebra maps}\label{InvSkewSymmBilMaps}

Except as noted in Lemmas~\ref{eigenvectors} and~\ref{le:orbitphixyphig}, Proposition~\ref{prop:directsum},
and Corollary~\ref{prop:Extension}, for the rest of~the paper,
let $G=S_n$ be the symmetric group with $n\geq 3$, let
$W \cong \mathbb{C}^n$ denote its natural permutation representation, and
consider the doubled permutation representation of $S_n$ on~$V=W^* \oplus W\cong \mathbb{C}^{2n}$.
Let $\mathcal{B}_y=\{y_1,\dots,y_n\}$ be the standard basis for $W$ and
$\mathcal{B}_x=\{x_1,\dots,x_n\}$ be the corresponding dual basis for $W^*$. Then the action
of $\sigma \in S_n$ is given by $\actby{\sigma}y_i=y_{\sigma(i)}$ and~$\actby{\sigma}x_i=x_{\sigma(i)}$. Recall that $W^* \cong \mathfrak{h}^* \oplus \iota^*$,
where $S_n$ acts trivially on the $1$-dimensional subspace $\iota^*$ of $W^*$
spanned by $x_{[n]}=\sum\limits_{i=1}^n x_i$ and by the standard reflection
representation on its $(n-1)$-dimensional orthogonal complement $\mathfrak{h}^*$,
and similarly $W \cong \mathfrak{h} \oplus \iota$.
In Remark~\ref{rmk:subrep1} and Section~\ref{DoubledStandardDefs} we also
consider the doubled standard representation of $S_n$ on the
subspace $\mathfrak{h}^* \oplus \mathfrak{h}$ spanned by{\samepage
\begin{gather}
\label{xbar}
\bigg\{\bar{x}_i:=x_i-\frac{1}{n}x_{[n]},\,\bar{y}_i:=y_i-\frac{1}{n}y_{[n]} \,\bigg|\, 1 \leq i \leq n\bigg\}
\end{gather}
or by $\{x_i-x_j,\, y_i-y_j \mid 1 \leq i,j \leq n\}$.}

In this section we identify all pre-Drinfeld orbifold algebra maps for $S_n$ with $n \geq 3$ acting by the doubled permutation representation on $W^*\oplus W$.
That is, we find all linear $2$-cochains~$\kappa^L$
satisfying the image condition~\eqref{PBWconditions-0}, the $G$-invariance condition~\eqref{PBWconditions-i}, and the mixed Jacobi identity $\psi\big(\kappa^L\big)=0$~\eqref{PBWconditions-ii}. To organize computations we make use of Lemma~\ref{ineedaname} relating Hochschild cohomology
and pre-Drinfeld orbifold algebra maps.

By Remark~\ref{remark:codim}, we need only consider group elements whose fixed point space has codimension zero or two.
Thus for $S_n$ acting by the doubled permutation representation we consider two cases: $\kappa^L$ supported only on the identity
and $\kappa^L$ supported only on the set of transpositions (which act as reflections on $W$ and bireflections on $W^*\oplus W$).

\subsection{Pre-Drinfeld orbifold algebra maps supported only on the identity}
We first prove a lemma that describes all $G$-invariant maps $\kappa_1\colon \bigwedge^2 V\to V\oplus \mathbb{C}$,
where $G$ is any finite group and $V$ is a permutation representation of $G$. Since properties~\eqref{PBWconditions-0}
and~\eqref{PBWconditions-ii}
are trivially satisfied when $g=1$, this will produce pre-Drinfeld orbifold algebra
maps supported only on the identity.
Recall that $\kappa_1$ is $G$-invariant if and only if
\begin{gather}
\label{G-inv}
\kappa_1(\actby{g}u,\actby{g}v)=\actbylessspace{g}(\kappa_1(u,v))
\end{gather}
for all $g$ in~$G$ and $u,v\in V$.
The following lemma shows that how $G$ acts on a set of representative basis
vector pairs determines a $G$-invariant linear cochain.

\begin{Lemma}
\label{eigenvectors}
	Suppose $G$ is a finite group acting on a complex vector space $V$ by
	a permutation representation. If~$\kappa^L_1$ is $G$-invariant, then the following two conditions hold for all $g$ in~$G$ and all basis vector pairs $v_i$ and $v_j$.
	\begin{enumerate}[label=$(\roman*)$]\itemsep=0pt
		\item
If $g$ swaps $v_i$ and $v_j$, then $\kappa^L_1(v_i,v_j)$ is an eigenvector of $g$ with eigenvalue $-1$.
		\item
If $g$ fixes $v_i$ and $v_j$, then the vector $\kappa^L_1(v_i,v_j)$ is in the fixed space $V^g$.
	\end{enumerate}
Suppose every ordered pair $v$, $w$ of basis vectors can be related to some $v_i$, $v_j$ in a set $S$ of~repre\-sentative basis vector pairs by $v=\actby{g}v_i$ and $w=\actby{g}v_j$ or $v=\actby{g}v_j$ and $w=\actby{g}v_i$ for some $g \in G$. If~$\kappa^L_1\colon S \to V$ satisfies~$(i)$ and~$(ii)$ for all representative pairs in $S$, then there is a unique way to extend $\kappa^L_1$ to be $G$-invariant on $\bigwedge^2 V$.
\end{Lemma}

\begin{proof}
Assume $\kappa^L_1$ is $G$-invariant. By~\eqref{G-inv}, if $g$ fixes both $v_i$ and $v_j$, then $\kappa^L_1(v_i,v_j)$ must be an element of $V^g$. And if $g$ swaps $v_i$ and $v_j$, then due to skew-symmetry, $\kappa^L_1(v_i,v_j)$ must
be a~$(-1)$-eigenvector of $g$.

Suppose $\kappa^L_1\colon \bigwedge^2 V\to V$ satisfies~($i$) and~($ii$) for a set of representative basis vector pairs. If~$v=\actby{g}v_i=\actby{h}v_i$ and $w=\actby{g}v_j=\actby{h}v_j$ for some representative pair $v_i$, $v_j$ and some $g, h \in G$, then $h^{-1}g$ fixes the basis vector pair $v_i$, $v_j$. Hence by~($ii$),
$\kappa^L_1(v_i,v_j)$ is in $V^{h^{-1}g}$ and $\actby{g}\kappa^L_1(v_i,v_j)=\actby{h}\kappa^L_1(v_i,v_j)$. If~instead $v=\actby{g}v_i=\actby{h}v_j$ and $w=\actby{g}v_j=\actby{h}v_i$, then $h^{-1}g$ swaps $v_i$ and $v_j$. Hence by~($i$),
$\actby{g}\kappa^L_1(v_i,v_j)=\actby{h}\kappa^L_1(v_j,v_i)$.
These imply that the unique way to extend $\kappa^L_1$ to be $G$-invariant is well-defined.
\end{proof}

We now apply this to the doubled permutation representation of $S_n$
on $W^*\oplus W$ equipped with the basis
$\mathcal{B}_x \cup \mathcal{B}_y$,
where $\mathcal{B}_y=\{y_1,\dots,y_n\}$ is the standard basis for $W$ and
$\mathcal{B}_x=\{x_1,\dots,x_n\}$ is the corresponding dual basis for $W^*$.
The following proposition summarizes the definitions
of~all $G$-invariant skew-symmetric bilinear maps, i.e., describes $\big(H_1^{2,0}\oplus H_1^{2,1}\big)^G$.

\begin{Proposition}
\label{prop:kappa1}
Let $S_n$ $(n \geq 3)$ act by the doubled permutation representation on $V=W^* \oplus W \cong \mathbb{C}^{2n}$ equipped with basis $\mathcal{B}_x \cup \mathcal{B}_y$.
The $S_n$-invariant linear and constant $2$-cochains
$\kappa^L_1\colon \bigwedge^2 V \to V$ and
$\kappa^C_1\colon \bigwedge^2 V\to \mathbb{C}$ are as given in Definition~$\ref{def:kappa1}$ in
terms of complex parameters $a_k, b_k$ for $1 \leq k \leq 7$ and $\alpha, \beta$ in $\mathbb{C}$ respectively.
\end{Proposition}

\begin{proof}
\textit{Linear cochains.}
If $V\cong \mathbb{C}^{2n}$ is the doubled permutation representation of $S_n$
and $\kappa^L_1\colon \bigwedge^2 V\to V$ is an $S_n$-invariant map,
then the value on any pair of basis vectors in $\mathcal{B}_x \cup \mathcal{B}_y=\{x_1,\dots,x_n,y_1,\dots,y_n\}$
can be obtained by acting by an appropriate permutation on one of the representative values
$\kappa^L_1(x_1,x_2)$, $\kappa^L_1(y_1,y_2)$, $\kappa^L_1(x_1,y_1)$, or $\kappa^L_1(x_1,y_2)$.

Consider $\kappa^L_1(x_1,x_2)$. The permutation $\sigma=(12)$ swaps $x_1$ and $x_2$, so by Lemma~\ref{eigenvectors}, $\kappa^L_1(x_1,x_2)$ must be a $(-1)$-eigenvector of $(12)$, i.e.,
a linear combination of the vectors $x_1-x_2$ and $y_1-y_2$. Both of these vectors are fixed by the group $S_{\{3,\ldots,n\}}$ of permutations that fix
both $x_1$ and $x_2$. Thus, for a choice of complex parameters $a_1$ and $b_1$,
we let
\begin{gather*}
\kappa^L_1(x_1,x_2)=a_1(x_1-x_2)+b_1(y_1-y_2).
\end{gather*}
A similar argument shows we can let
\begin{gather*}
\kappa^L_1(y_1,y_2)=a_2(x_1-x_2)+b_2(y_1-y_2)
\end{gather*}
for some choice of complex parameters $a_2$ and $b_2$.

Consider $\kappa^L_1(x_1,y_1)$. There are no permutations that will swap $x_1$ and $y_1$. The group of~per\-mu\-tations that fix both $x_1$ and $y_1$ is $S_{\{2,\ldots,n\}}$, so by Lemma~\ref{eigenvectors},
$\kappa^L_1(x_1,y_1)$ must be an element of the subspace
\begin{gather*}
V^{S_{\{2,\ldots,n\}}}=\Span\big\{x_1, x_{[n]}, y_1, y_{[n]}\big\}.
\end{gather*}
We define $\kappa^L_1(x_1,y_1)$ to be a linear combination of the basis elements, using complex parameters $a_3$, $a_4$, $b_3$, and $b_4$ as weights.
Orbiting yields the definition
\begin{gather*}
\kappa^L_1(x_i,y_i)=a_3x_i+a_4x_{[n]}+b_3y_i+b_4y_{[n]}
\end{gather*}
for $1\leq i\leq n$.

Consider $\kappa^L_1(x_1,y_2)$. There are no permutations that will swap $x_1$ and $y_2$. The group of~per\-mu\-tations that fix both $x_1$ and $y_2$ is $S_{\{3,\ldots,n\}}$, so once again by Lemma~\ref{eigenvectors}, $\kappa^L_1(x_1,y_2)$ must be an element of
the subspace
\begin{gather*}
V^{S_{\{3,\ldots,n\}}}=\Span\big\{x_1,x_2, x_{[n]}, y_1, y_2, y_{[n]}\big\}.
\end{gather*}
We define $\kappa^L_1(x_1,y_2)$ to be a linear combination of the basis elements using complex parameters $a_5$, $a_6$, $a_7$, $b_5$, $b_6$, and $b_7$ as weights.
Orbiting yields the definition
\begin{gather*}
\kappa^L_1(x_i,y_j)=a_5x_i+a_6x_j+a_7x_{[n]}+b_5y_i+b_6y_j+b_7y_{[n]}
\end{gather*}
for $1\leq i,j\leq n$ with $i\neq j$.

\textit{Constant cochains.}
By comparison there is only a two-parameter family of $G$-invariant constant cochains.
First, using~\eqref{G-inv}, if a constant cochain $\kappa_1^C\colon \bigwedge^2 V\to \mathbb{C}$ is
$S_n$-invariant and some element $g \in S_n$ swaps $v_i$ and $v_j$, then $\kappa_1^C(v_i,v_j)=0$.
Thus $\kappa_1^C(x_i,x_j)=\kappa_1^C(y_i,y_j)=0$.
Also due to $S_n$-invariance, the value of $\kappa_1^C(x_i,y_j)$ only depends on whether $i=j$ or $i\neq j$, so for $\alpha,\beta\in\mathbb{C}$, we can let $\kappa_1^C(x_i,y_i)=\alpha$ and $\kappa_1^C(x_i,y_j)=\beta$ for $i\neq j$. This shows we have a $2$-dimensional space of $S_n$-invariant maps $\kappa_1^C$.
\end{proof}

\begin{Remark}
It is also possible to confirm the dimensions for the linear and constant invariant cochains by using the equivalences from Lemma~\ref{ineedaname} and calculating inner products of characters. If~$n\geq 3$, then, omitting details, we find
for the linear cochains,
\begin{gather*}
\dim\big(H_1^{2,1}\big)^{S_n}=\dim\Big(V\otimes\bigwedge\nolimits^{\! 2} V^*\Big)^{S_n}=\big\langle \chi_{\iota},\chi_{_V}\chi_{\bigwedge^2 V}\big\rangle=\big\langle\chi_{_V},\chi_{\bigwedge^2 V}\big\rangle =14,
\end{gather*}
and for the constant cochains,
\begin{gather*}
\dim \big(H_1^{2,0}\big)^{S_n}=\dim\Big(\bigwedge\nolimits^{\! 2} V^*\Big)^{S_n}=\langle\chi_{\iota},\chi_{\bigwedge^2 V} \rangle=2,
\end{gather*}
as expected.
\end{Remark}

\subsection{Pre-Drinfeld orbifold algebra maps supported only off the identity}

The following proposition describes $\big(H_g^{2,0}\oplus H_g^{2,1}\big)^G$
where $g$ is a transposition.

\begin{Proposition}
\label{prop:kapparefl}
Let $S_n$ $(n \geq 3)$ act by the doubled permutation representation on $V=W^* \oplus W \cong \mathbb{C}^{2n}$, equipped with the
basis $\mathcal{B}_x \cup \mathcal{B}_y$.
The $S_n$-invariant linear and constant $2$-cochains that sastify the mixed Jacobi identity and are supported only on transpositions are the maps of the form given in Definition~$\ref{def:kapparefl}$.
\end{Proposition}

\begin{proof}
We find the centralizer invariants $\big(H_g^{2,0}\oplus H_g^{2,1}\big)^{Z(g)}$ when $g$ is a transposition.
Let $g=(12)$ and first note that the centralizer of $g$ is $Z(g)=\langle (12) \rangle \times S_{\{3,\dots,n\}}$, the fixed point space of $g$ is
\begin{gather*}
V^g=\Span\{x_1+x_2,x_3,\dots,x_n,y_1+y_2,y_3,\dots,y_n\},
\end{gather*}
and the orthogonal complement is
\begin{gather*}
(V^g)^\perp=\Span\{x_1-x_2,y_1-y_2\}.
\end{gather*}
A basis for $\bigwedge^2 ((V^g)^{\perp})^*$
is the volume form
	\begin{gather*}
	\vol{g}:=(x_1^*-x_2^*)\wedge (y_1^*-y_2^*).
	\end{gather*}
Note that $Z(g)$ acts trivially on $\vol{g}$, so
 	\begin{gather*}
	\big(H_g^{2,0}\big)^{Z(g)}=H_g^{2,0}=\bigwedge\nolimits^{\! 2} \big((V^g)^{\perp}\big)^* \otimes \mathbb{C} g=\Span\big\{\vol{g} \otimes (12) \big\}
	\end{gather*}
and
\begin{align*}
\big(H_g^{2,1}\big)^{Z(g)}& =V^{Z(g)}\otimes\bigwedge\nolimits^{\! 2} \big((V^g)^{\perp}\big)^* \otimes \mathbb{C} g
\\
&=\Span\bigg\{v \otimes \vol{g} \otimes (12) \,\bigg\vert\, v \in \bigg\{x_1+x_2, \sum\limits_{i=3}^n x_i, y_1+y_2, \sum\limits_{i=3}^n y_i \bigg\}\bigg\}.
\end{align*}

After orbiting the centralizer invariants to obtain~$G$-invariants (see the end of Section~4.1 in~\cite{FGK} for more detail), these yield the description in Definition~\ref{def:kapparefl} for the cochain $\kappa_{\refl}=\sum_{(ij) \in S_n} \big(\kappa^C_{(ij)}+\kappa^L_{(ij)}\big)\otimes (ij)$ supported off the identity.
\end{proof}

\subsection{Pre-Drinfeld orbifold algebra maps}
By Lemma~\ref{ineedaname} and Remark~\ref{remark:codim}, the polynomial
degree one elements of Hochschild $2$-cohomology found in
Propositions~\ref{prop:kappa1} and~\ref{prop:kapparefl}
provide a description of all pre-Drinfeld orbifold algebra maps.

\begin{Corollary}
\label{preDOAmaps}
The pre-Drinfeld orbifold algebra maps for $S_n$ $(n\geq 3)$ acting
by the doubled permutation representation on
$V=W^* \oplus W \cong \mathbb{C}^{2n}$ are the linear $2$-cochains
$\kappa^L=\kappa^L_1+\kappa^L_{\refl}$ for $\kappa^L_1$ described in
terms of the parameters $a_1,\dots,a_7$, $b_1,\dots,b_7$ as in
Definition~$\ref{def:kappa1}$ and
$\kappa^L_{\refl}$ controlled by the
parameters $a$, $a^\perp$, $b$, $b^\perp$
as in Definition~$\ref{def:kapparefl}$.
\end{Corollary}

In Theorems~\ref{thm:LOAMapsDoubledPerm} and~\ref{thm:DOAMapsDoubledPerm}
we will characterize when the maps $\kappa^L_1$ and $\kappa^L_{\refl}$ lift separately
to Drinfeld orbifold algebra maps and in Theorem~\ref{thm:DOAMapsCombined} we will
show it is also possible to lift $\kappa^L_1+\kappa^L_{\refl}$. Any two lifts
of a particular pre-Drinfeld orbifold algebra map must differ by a~con\-stant
$2$-cochain that satisfies the mixed Jacobi identity. Lemma~\ref{ineedaname} and the results in this section yield the
following corollary describing these maps.

\begin{Corollary}\label{mixedJacobi}
For $S_n$ $(n\geq 3)$ acting on $V=W^* \oplus W \cong \mathbb{C}^{2n}$ by the
doubled permutation representation, the
$S_n$-invariant constant $2$-cochains satisfying the mixed Jacobi identity
are the maps $\kappa^C=\kappa^C_1+\kappa^C_{\refl}$ with $\kappa^C_1$ given in
terms of parameters $\alpha$ and $\beta$ in Definition~$\ref{def:kappa1}$
and $\kappa^C_{\refl}$ described using parameter $c$ in Definition~$\ref{def:kapparefl}$.
\end{Corollary}

\subsection{Definitions of linear and constant cochains}
For convenience we collect here the definitions of the components of the maps
determined in~Pro\-positions~\ref{prop:kappa1} and~\ref{prop:kapparefl} and that
will be needed to lift $\kappa^L_1$ in~Section~\ref{InvLieBrackets} and
$\kappa^L_{\refl}$ in~Section~\ref{Lifting}.

Some parts of the descriptions below involve sums of basis vectors
over subsets of $[n]=\{1,\dots,n\}$. For~$I \subseteq [n]$ let $v_I=\sum_{i \in I} v_i$,
where $v$ stands for $x$ or $y$ and at times we omit the set braces in $I$.
Let $v_I^\perp$ denote the complementary vector $v_{[n]}-v_I$.
In all three definitions, $S_n$ $(n \geq 3)$ acts by the doubled permutation representation on $V=W^* \oplus W \cong \mathbb{C}^{2n}$ equipped with basis $\mathcal{B}_x \cup \mathcal{B}_y=\{x_1,\dots,x_n,y_1,\dots,y_n\}$.

\begin{Definition}[cochains supported only on the identity]\label{def:kappa1}
Given complex parameters $a_k, b_k$ for $1 \leq k \leq 7$ and $\alpha, \beta$ in $\mathbb{C}$,
let $\kappa^L_1\colon \bigwedge^2 V \to V$ and
$\kappa^C_1\colon \bigwedge^2 V\to \mathbb{C}$
be the $S_n$-invariant maps defined by
\begin{gather}
\label{kappaLxx}
\kappa^L_1(x_i,x_j)=a_1(x_i-x_j)+b_1(y_i-y_j), \\
\label{kappaLyy}
\kappa^L_1(y_i,y_j)=a_2(x_i-x_j)+b_2(y_i-y_j), \\
\label{kappaLxiyi}
\kappa^L_1(x_i,y_i)=a_3x_i+a_4\n{x}+b_3y_i+b_4\n{y}, \\
\label{kappaLxiyj}
\kappa^L_1(x_i,y_j)=a_5x_i+a_6x_j+a_7\n{x}+b_5y_i+b_6y_j+b_7\n{y},
\end{gather}
and
\begin{gather*}
\kappa^C_1(x_i,x_j)=\kappa^C_1(y_i,y_j)=0, \qquad \kappa^C_1(x_i,y_i)=\alpha, \qquad \kappa^C_1(x_i,y_j)=\beta,
\end{gather*}
where $1 \leq i \neq j \leq n$.
\end{Definition}

\begin{Definition}[cochains supported only on transpositions]\label{def:kapparefl}
Let $a$, $a^\perp$, $b$, $b^\perp$, $c$ be complex parameters and let $T$ be the set of transpositions in~$S_n$.
Define a linear $2$-cocycle $\kappa^L_{\refl}=\sum_{g\in T}\kappa_g^Lg$, where for $g=(rs)$, the component $\kappa_g^L\colon \bigwedge^2 V\to V$ is defined for $1\leq i, j\leq n$ by
\begin{gather*}
\kappa_g^L(x_i,x_j)=\kappa_g^L(y_i,y_j)=0
\end{gather*}
and
\begin{gather*}
\kappa_g^L(x_i,y_j)=\begin{cases}
 ax_{r,s}+a^{\perp}x_{r,s}^{\perp}+by_{r,s}+b^{\perp}y_{r,s} & \text{if}\quad i=j\quad \text{is in}\quad \{r,s\},
 \\
 -(ax_{r,s}+a^{\perp}x_{r,s}^{\perp}+by_{r,s}+b^{\perp}y_{r,s}) & \text{if}\quad \{i,j\}=\{r,s\},
 \\
 0 & \text{otherwise.}
\end{cases}
\end{gather*}
Similarly, the $g=(rs)$ component of the constant $2$-cocycle $\kappa^C_{\refl}=\sum_{g\in T}\kappa_g^Cg$ is defined for $1\leq i, j\leq n$ by
\begin{gather*}
\kappa_g^C(x_i,x_j)=\kappa_g^C(y_i,y_j)=0
\end{gather*}
and
\begin{gather*}
\kappa_g^C(x_i,y_j)=\begin{cases}
c & \text{if}\quad i=j\quad \text{is in}\quad \{r,s\},
\\
-c & \text{if}\quad \{i,j\}=\{r,s\},
 \\
0 & \text{otherwise}.
\end{cases}
\end{gather*}
\end{Definition}

Lastly, we define a constant $2$-cochain $\kappa^C_{\tri}$ which we use to lift $\kappa^L_{\refl}$ in~Section~\ref{ObstructionDetails}. The map $\kappa^C_{\tri}$ is not a Hochschild $2$-cocycle but rather is
based on the form of $\phi(\kappa^L_{\refl},\kappa^L_{\refl})$ in Propositions~\ref{prop:zerocases} and~\ref{prop:3-cycle} and is constructed to ensure $\phi(\kappa^L_{\refl},\kappa^L_{\refl})=2\psi(\kappa^C_{\tri})$ as in
Proposition~\ref{prop:phiL2L2=2psiC3}, thereby
clearing the first obstruction to lifting $\kappa^L_{\refl}$. The cochain $\kappa^C_{\tri}$ will also clear the second obstruction to lifting $\kappa^L_{\refl}$, as verified in Lemma~\ref{lemma:simplificationC3L2}.

\begin{Definition}[cochains supported only on $3$-cycles]\label{def:kappatri}
Define an $S_n$-invariant map $\kappa^C_{\tri}=\sum_{g\in S_n}\kappa^C_gg$ with component maps $\kappa^C_{g}\colon \bigwedge^2 V\to \mathbb{C}$. If~$g$ is not a $3$-cycle, let $\kappa^C_g\equiv 0$. For~a~$3$-cycle $g=(i~j~k)$, define the outcome of $\kappa_{g}^C$ on a pair of basis vectors to be zero unless the indices are two distinct elements of $\{i,j,k\}$, in which case the outcome is defined by the following (and skew-symmetry):
\begin{gather*}
\kappa_{g}^C(x_i,x_j)=\kappa_{g}^C(x_j,x_k)=\kappa_{g}^C(x_k,x_i)=\big(b^{\perp}-b\big)^2
\end{gather*}
and
\begin{gather*}
\kappa_{g}^C(y_i,y_j)=\kappa_{g}^C(y_j,y_k)=\kappa_{g}^C(y_k,y_i)=\big(a^{\perp}-a\big)^2
\end{gather*}
and
\begin{gather*}
\kappa_{g}^C(x_i,y_j)=\kappa_{g}^C(y_j,x_k)=\kappa_{g}(x_k,y_i)
=\kappa_{g}^C(y_i,x_j)=\kappa_{g}^C(x_j,y_k)=\kappa_{g}(y_k,x_i)
\\ \hphantom{\kappa_{g}^C(x_i,y_j)}
{}=-\big(a^{\perp}-a\big)\big(b^{\perp}-b\big).
\end{gather*}
\end{Definition}

\section[Lie orbifold algebra maps that deform S(W* + W)Sn]
{Lie orbifold algebra maps that deform $\boldsymbol{S(W^* \oplus W)\# S_n}$}\label{InvLieBrackets}

In Section~\ref{InvSkewSymmBilMaps}, as summarized in Proposition~\ref{prop:kappa1} and Definition~\ref{def:kappa1},
we determined the
pre-Drinfeld orbifold algebra maps $\kappa^L_1$ supported only on the identity.
Here we find conditions under which these maps also endow
$V$ with a Lie algebra structure~--- i.e., under which they lift to Lie orbifold
algebra maps because there exists a constant $2$-cochain $\kappa^C$
such that $\kappa=\kappa^L_1+\kappa^C$ also satisfies the remaining
properties~\eqref{PBWconditions-iii} and~\eqref{PBWconditions-iv}.

Our main goal is to write down conditions on the parameters involved in the definitions
of $\kappa^L_1$, $\kappa^C_1$, and $\kappa^C_{\refl}$ such that properties~\eqref{PBWconditions-iii}
and~\eqref{PBWconditions-iv} hold, or in other words, such that
$\phi\big(\kappa^L_1,\kappa^L_1\big)=2\psi\big(\kappa^C_1+\kappa^C_{\refl}\big)$ and
$\phi\big(\kappa^C_1+\kappa^C_{\refl},\kappa^L_1\big)=0$.
Since $2\psi\big(\kappa^C_1+\kappa^C_{\refl}\big)=0$, we have that
$\kappa^C_1+\kappa^C_{\refl}$ clears both the first and second obstructions and
the map $\kappa^L_1$ gives rise to a Lie orbifold algebra if and only if
$\phi\big(\kappa^L_1,\kappa^L_1\big)=\phi\big(\kappa^C_1+\kappa^C_{\refl},\kappa^L_1\big)=0$. We~use this to arrive at characterizing PBW conditions on parameters as summarized in
the proof of
Theorem~\ref{thm:LOAMapsDoubledPerm}, which states that $\kappa^L_1$ can
be lifted to $\kappa=\kappa^L_1+\kappa^C_1+\kappa^C_{\refl}$ precisely when a list of $22$
homogeneous quadratic conditions in $17$ parameters hold.

It will be convenient along the way to also consider $\phi\big(\kappa^L_{\refl},\kappa^L_1\big)$,
for use in Theorem~\ref{thm:DOAMapsCombined}, by using $*$ to denote
either $C$ or $L$ and $x$ to denote either a transposition or the identity and
computing, for $v_1,v_2,v_3 \in V$,
\begin{gather*}
\phi_{x,1}^*(v_1,v_2,v_3):=\kappa^*_x\big(v_1,\kappa^L_1(v_2,v_3)\big)
+\kappa^*_x\big(v_2,\kappa^L_1(v_3,v_1)\big)+\kappa^*_x\big(v_3,\kappa^L_1(v_1,v_2)\big)
\end{gather*}
as uniformly as possible.
This notation omits a factor of two (and hence differs from that in~\cite{FGK})
because $\psi(\kappa^C_1+\kappa^C_{\refl})=0$ means the factor of 2 is irrelevant to
clearing the first obstruction and it is also irrelevant to clearing
the second obstruction.

First note that due to bilinearity and skew-symmetry it suffices to compute $\phi_{x,1}^*$,
with $x$ equal to the identity or a transposition, on basis triples of six main types
for $1 \leq i, j, k \leq n$, where $n\geq 3$.
\begin{alignat*}{3}
&\text{1.\ All basis vectors in $W$ or in $W^*$ and $i$, $j$, $k$ distinct:} &\hspace{10pt}
& (x_i,x_j,x_k), & &\quad (y_i,y_j,y_k). \\
&\text{2.\ Two basis vectors in $W$ or in $W^*$ and $i$, $j$, $k$ distinct:} &\hspace{10pt}
& (x_i,x_j,y_k), & &\quad (y_i,y_j,x_k). \\
&\text{3.\ Two basis vectors in $W$ or $W^*$ and $i$, $j$ distinct:} &\hspace{10pt}
& (x_i,x_j,y_j), & &\quad (y_i,y_j,x_j).
\end{alignat*}
This is done in the next three subsections.

\subsection[All basis vectors in W or in W* and three distinct indices]
{All basis vectors in $\boldsymbol W$ or in $\boldsymbol{W^*}$ and three distinct indices}\label{AllVecInW}

For any distinct indices $i$, $j$, $k$ with $1 \leq i, j, k \leq n$, we have
\begin{gather*}
\phi_{x,1}^*(x_i,x_j,x_k)=\kappa^*_x\big(x_i,\kappa^L_1(x_j,x_k)\big)+\kappa^*_x\big(x_j,\kappa^L_1(x_k,x_i)\big) +\kappa^*_x\big(x_k,\kappa^L_1(x_i,x_j)\big).
\end{gather*}
Using bilinearity, skew-symmetry, and Definitions~\ref{def:kappa1}
and~\ref{def:kapparefl} of $\kappa_1$ and
$\kappa_{\refl}$ yields for $x$ either the identity or
any transposition,
\begin{gather*}
\kappa^*_x(x_i,x_j-x_k)+\kappa^*_x(x_j,x_k-x_i)+\kappa^*_x(x_k,x_i-x_j)
\\ \qquad
{}=2[\kappa^*_x(x_i,x_j)+\kappa^*_x(x_j,x_k)+\kappa^*_x(x_k,x_i)]=0,	
\end{gather*}
and
\begin{gather*}
\kappa^*_x(x_i,y_j-y_k)+\kappa^*_x(x_j,y_k-y_i)+\kappa^*_x(x_k,y_i-y_j)=0.
\end{gather*}
Combining these shows that $\phi_{x,1}^*(x_i,x_j,x_k)=0$
and similarly $\phi_{x,1}^*(y_i,y_j,y_k)=0$, for any (distinct $i$, $j$, $k$ with) $1 \leq i, j, k \leq n$,
for $x$ either the identity or a transposition, and with $*=C$ or~$*=L$.
Thus this case imposes no restrictions on any parameters.

\subsection[Two basis vectors in W or W* and three distinct indices]
{Two basis vectors in $\boldsymbol W$ or $\boldsymbol{W^*}$ and three distinct indices}

For any distinct indices $i$, $j$, $k$ with $1 \leq i, j, k \leq n$,
using the definition of $\kappa^L_1$, bilinearity, and skew-symmetry yields
\begin{gather*}
\phi_{x,1}^*(x_i,x_j,y_k)
=2a_5 \kappa^*_x(x_i,x_j) + a_6 \kappa^*_x(x_i-x_j,x_k) + a_7 \kappa^*_x\big(x_i-x_j,x_{[n]}\big)
\\ \phantom{\phi_{x,1}^*(x_i,x_j,y_k)=}
{}-b_1 \kappa^*_x(y_i-y_j,y_k)+b_5\left(\kappa^*_x(x_i,y_j) - \kappa^*_x(x_j,y_i)\right)
\\ \phantom{\phi_{x,1}^*(x_i,x_j,y_k)=}
{}+(b_6 - a_1)\kappa^*_x(x_i-x_j,y_k) + b_7 \kappa^*_x\big(x_i-x_j,y_{[n]}\big).
\end{gather*}

When $x=1$ and $*=C$, since $\kappa^C_1(v,w)=0$ when $v,w \in W$ or $v,w \in W^*$, we have
\begin{gather*}
\phi_{1,1}^C(x_i,x_j,y_k)=b_5(\beta-\beta)+(b_6-a_1)(\beta-\beta) +b_7\left(\alpha-\alpha+(n-1)(\beta-\beta)\right)=0,
\end{gather*}
and when $x=1$ and $*=L$ using the definition of $\kappa^L_1$ yields
\begin{gather*}
\phi_{1,1}^L(x_i,x_j,y_k)=\gamma_1(x_i-x_j)+\gamma_2(y_i-y_j),
\end{gather*}
where
\begin{gather*}
\gamma_1=a_1(a_5 + a_6 + na_7) - b_1 a_2 - b_5 a_6 + a_5 (b_5 + b_6 + nb_7) + b_7 (a_3 - a_5 - a_6),
\\
\gamma_2=b_1(a_5 + a_6 + na_7) - b_1 b_2 + b_1 a_5 + b_5 (b_5 - a_1 + nb_7) + b_7 (b_3 - b_5 - b_6).
\end{gather*}
When $x=g$ is a transposition, by Definition~\ref{def:kapparefl} we have
that
\begin{align*}
\phi_{g,1}^*(x_i,x_j,y_k)&=(b_6 - a_1)\kappa^*_g(x_i-x_j,y_k) \\
&=
\begin{cases}
 \pm(b_6 - a_1)\kappa^*_g(x_l,y_k), & \text{if\ \ $g=(lk)$\ \ with\ \ $l=i$\ \ or\ \ $l=j$\ \ respectively}, \\
 0, & \text{otherwise},
\end{cases}
\end{align*}
and we define
\begin{gather*}
\gamma_3=-(b_6 - a_1)\kappa_{(jk)}^C(x_j,y_k)=c(b_6 - a_1).
\end{gather*}

Interchanging the roles of $x$ and $y$ and recomputing yields that for
distinct $i$, $j$, $k$ with $1 \leq i, j, k \leq n$,
\begin{gather*}
\phi_{x,1}^*(y_i,y_j,x_k)
=-2b_6 \kappa^*_x(y_i,y_j) - b_5 \kappa^*_x(y_i-y_j,y_k) - b_7 \kappa^*_x\big(y_i-y_j,y_{[n]}\big)
\\ \hphantom{\phi_{x,1}^*(y_i,y_j,x_k)=}
{}- a_2 \kappa^*_x(x_i-x_j,x_k)-a_6\big(\kappa^*_x(x_i,y_j) - \kappa^*_x(x_j,y_i)\big)
\\ \hphantom{\phi_{x,1}^*(y_i,y_j,x_k)=}
{}+(b_2 + a_5)\kappa^*_x(x_k,y_i-y_j) + a_7 \kappa^*_x\big(x_{[n]},y_i-y_j\big).
\end{gather*}

In particular,
\begin{gather*}
\phi_{1,1}^C(y_i,y_j,x_k)
=0 \qquad \text{and} \qquad \phi_{1,1}^L(y_i,y_j,x_k)
=\gamma_4(x_i-x_j)+\gamma_5(y_i-y_j),
\end{gather*}
where
\begin{gather*}
\gamma_4=-a_2(b_6 + b_5 + nb_7) - a_2 a_1 - a_2 b_6 + a_6(a_6 + b_2 + na_7) + a_7(a_3 - a_6 - a_5), \\
\gamma_5=-b_2(b_6 + b_5 + nb_7) - a_2 b_1 - a_6 b_5 + b_6(a_6 + a_5 + na_7) + a_7(b_3 - b_6 - b_5).
\end{gather*}
When $x=g$ is a transposition, we have
\begin{align*}
\phi_{g,1}^*(y_i,y_j,x_k)&=(b_2 + a_5)\kappa^*_g(x_k,y_i-y_j) \\
&=
\begin{cases}
 \pm(b_2 + a_5)\kappa^*_g(x_k,y_l) & \text{if\ \ $g=(lk)$\ \ with\ \ $l=i$\ \ or\ \ $l=j$\ \ respectively}, \\
 0 & \text{otherwise},
\end{cases}
\end{align*}
and we define
\begin{gather*}
\gamma_6=-(b_2 + a_5)\kappa_g^C(x_k,y_j)=c(b_2 + a_5).
\end{gather*}

\subsection[Two basis vectors in W or W^* and two distinct indices]
{Two basis vectors in $\boldsymbol W$ or $\boldsymbol {W^*}$ and two distinct indices}\label{TwoVecInWTwoIndices}

For any distinct indices $i$, $j$ with $1 \leq i, j \leq n$, we have
\begin{gather*}
\phi_{x,1}^*(x_i,x_j,y_j)
=(a_3+a_5) \kappa^*_x(x_i,x_j) - b_1 \kappa^*_x(y_i,y_j) + a_4 \kappa^*_x(x_i,x_{[n]}) - a_7 \kappa^*_x(x_j,x_{[n]})
\\ \hphantom{\phi_{x,1}^*(x_i,x_j,y_j)=}
{}+(-a_1+b_3) \kappa^*_x(x_i,y_j) - b_5 \kappa^*_x(x_j,y_i) + (a_1-b_6) \kappa^*_x(x_j,y_j)
\\ \hphantom{\phi_{x,1}^*(x_i,x_j,y_j)=}
{}+ b_4 \kappa^*_x(x_i,y_{[n]}) - b_7 \kappa^*_x(x_j,y_{[n]}).
\end{gather*}

In particular if $x=1$ and $*=C$ then we set
\begin{gather*}
\gamma_7=\phi_{1,1}^C(x_i,x_j,y_j)=\alpha(a_1-b_6+b_4-b_7)-\beta(a_1-b_3+b_5-(n-1)(b_4-b_7)),
\end{gather*}
and if $x=1$ and $*=L$ then
\begin{gather*}
\phi_{1,1}^L(x_i,x_j,y_j)
=\gamma_8 x_{[n]} +\gamma_9 x_i +\gamma_{10} x_j +\gamma_{11} y_{[n]} +\gamma_{12} y_i -\gamma_{13}y_j,
\end{gather*}
where
\begin{gather*}
\gamma_8=a_7(b_3 - b_5) - a_4 b_6 + (b_4-b_7)(a_4 + (n-1)a_7 + a_6), \\
\gamma_9=a_1(a_3 + na_4) + a_5(b_3 + nb_4) + b_4(a_3 - a_5 - a_6) - a_2 b_1 - b_5 a_6, \\
\gamma_{10}=-a_1(a_5 + a_6 + na_7) - a_5(b_5 + nb_7) - b_7(a_3 - a_5 - a_6) + a_2 b_1 + b_3 a_6 - a_3 b_6, \\
\gamma_{11}=b_7(b_3 - b_5) - b_4 b_6 + (b_4 - b_7)(a_1 + b_4 + (n-1) b_7 + b_6) - b_1(a_4 - a_7), \\
\gamma_{12}=b_1(a_5 - b_2 + a_3 + na_4) - b_5(a_1 + b_6 - b_3 - nb_4) + b_4(b_3 - b_5 - b_6),
\\
\gamma_{13}=b_1(-b_2 + a_5 + a_3 + na_7) + b_5(b_5 + nb_7) + b_7(b_3 - b_5 - b_6) - a_1(b_3 - b_6).
\end{gather*}
For $x=g$ a transposition, we have that
\begin{align*}
\phi_{g,1}^*(x_i,x_j,y_j)&=(-a_1 + b_3)\kappa^*_g(x_i,y_j) - b_5\kappa^*_g(x_j,y_i) + (a_1 - b_6)\kappa^*_g(x_j,y_j)
\\
&=\begin{cases}
 (2a_1 - b_3 + b_5 - b_6)\kappa^*_g(x_j,y_j), & \text{if\quad $g=(ij)$,} \\
 0, & \text{otherwise},
\end{cases}
\end{align*}
and we define
\begin{gather*}
\gamma_{14}=(2a_1 - b_3 + b_5 - b_6)\kappa^C_g(x_j,y_j)=c(2a_1 - b_3 + b_5 - b_6).
\end{gather*}

Interchanging the roles of $x$ and $y$ and recomputing 
yields that for any distinct indices~$i$ and~$j$ with $1 \leq i, j \leq n$,
\begin{gather*}
\phi_{x,1}^*(y_i,y_j,x_j)
=-(b_3+b_6) \kappa^*_x(y_i,y_j) - a_2 \kappa^*_x(x_i,x_j) - b_4 \kappa^*_x(y_i,y_{[n]}) + b_7 \kappa^*_x(y_j,y_{[n]})
\\ \hphantom{\phi_{x,1}^*(y_i,y_j,x_j)=}
{}+(a_3+b_2) \kappa^*_x(x_j,y_i) - a_6 \kappa^*_x(x_i,y_j) - (a_5+b_2) \kappa^*_x(x_j,y_j)
\\ \hphantom{\phi_{x,1}^*(y_i,y_j,x_j)=}
{}+ a_4 \kappa^*_x(x_{[n]},y_i) - a_7 \kappa^*_x(x_{[n]},y_j).
\end{gather*}
In particular, we set
\begin{gather*}
\gamma_{15}=\phi_{1,1}^C(y_i,y_j,x_j)=\alpha(-b_2-a_5+a_4-a_7)+\beta(b_2+a_3-a_6+(n-1)(a_4-a_7)),
\end{gather*}
 and
\begin{align*}
\phi_{1,1}^L(y_i,y_j,x_j)
&=\gamma_{16} x_{[n]} +\gamma_{17} x_i +\gamma_{18} x_j +\gamma_{19}y_{[n]} +\gamma_{20}y_i +\gamma_{21}y_j,
\end{align*}
where
\begin{gather*}
\gamma_{16}=a_7(a_3 - a_6) - a_4 a_5 + (a_4-a_7)(-b_2 + a_4 + (n-1)a_7 + a_5) + a_2(b_4 - b_7), \\
\gamma_{17}=-a_2(a_1 + b_6 + b_3 + nb_4) + a_6(b_2 - a_5 + a_3 + na_4) + a_4(a_3 - a_5 - a_6), \\
\gamma_{18}=a_2(a_1 + b_6 + b_3 + nb_7) - a_6(a_6 + na_7) - a_7(a_3 - a_5 - a_6) - b_2(a_3 - a_5), \\
\gamma_{19}=b_7(a_3 - a_6) - b_4 a_5 + (a_4 - a_7)(b_4 + (n-1) b_7 + b_5), \\
\gamma_{20}=-b_2(b_3 + nb_4) + b_6(a_3 + na_4) + a_4(b_3 - b_5 - b_6) - a_2 b_1 - a_6 b_5, \\
\gamma_{21}=b_2(b_5 + b_6 + nb_7) - b_6(a_6 + na_7) - a_7(b_3 - b_5 - b_6) + a_2 b_1 + a_3 b_5 - b_3 a_5.
\end{gather*}
Lastly, when $x=g$ is a transposition, we have
\begin{align*}
\phi_{g,1}^*(x_i,x_j,y_j)&=(a_3 + b_2)\kappa^*_g(x_j,y_i) - a_6\kappa^*_g(x_i,y_j) - (a_5 + b_2)\kappa^*_g(x_j,y_j) \\
&=
\begin{cases}
 -(2b_2 + a_3 + a_5 - a_6)\kappa^*_g(x_j,y_j), & \text{if}\quad g=(ij), \\
 0, & \text{otherwise,}
\end{cases}
\end{align*}
and we define
\begin{gather*}
\gamma_{22}=(2b_2 + a_3 + a_5 - a_6)\kappa^*_g(x_j,y_j)=c(2b_2 + a_3 + a_5 - a_6).
\end{gather*}

\subsection{Lie orbifold algebra maps}

We now use the calculations in~Sections~\ref{AllVecInW}--\ref{TwoVecInWTwoIndices}
to describe all Drinfeld orbifold algebra maps with linear part
supported only on the identity, i.e., all Lie orbifold algebra maps
$\kappa^L_1+\kappa^C$. The~corresponding Lie orbifold algebras are
described in Theorem~\ref{Lieorbifold}.

\begin{Theorem}
\label{thm:LOAMapsDoubledPerm}
Let $S_n$ $(n \geq 3)$ act on $V=W^* \oplus W \cong \mathbb{C}^{2n}$ by the
doubled permutation representation, and let $\kappa^L_1$ and $\kappa^C_1$
be as described in Definition~$\ref{def:kappa1}$ and
$\kappa^C_{\refl}$ be as in Definition~$\ref{def:kapparefl}$ with
complex parameters $a_1,\dots,a_7$, $b_1,\dots,b_7$, $\alpha$, $\beta$, and $c$.
The Lie orbifold algebra maps are precisely the maps of the form
$\kappa=\kappa^L_1+\kappa^C_1+\kappa^C_{\refl}$ satisfying the
conditions $\gamma_i=0$ for $1 \leq i \leq 22$.
\end{Theorem}

\begin{proof}
Let $\kappa^L$ be a pre-Drinfeld orbifold algebra map supported only on the
identity. By~Corollary~\ref{preDOAmaps} we know $\kappa^L=\kappa^L_1$ is as given
in terms of $a_i$ and $b_i$ for $1 \leq i \leq 7$ in Definition~\ref{def:kappa1}.

In considering property~\eqref{PBWconditions-iii} when $g=1$, note that
$\psi_1 \equiv 0$ for any $\kappa^C$, so we must have $\phi_{1,1}\equiv 0$ as well.
The result of computing $\phi\big(\kappa^L_1,\kappa^L_1\big)$ is given
in~Sections~\ref{AllVecInW}--\ref{TwoVecInWTwoIndices} as the
values of the various $\phi^L_{1,1}(u,v,w)$.
These show that $\phi_{1,1}\equiv 0$
precisely when the parameters $a_i$, $b_i$ for $i=1,\dots,7$
satisfy the conditions
\begin{equation}
\label{ClearFirstObstr}
\gamma_1=\gamma_2=\gamma_4=\gamma_5=0\qquad \text{and}\qquad
\gamma_i=0 \quad \text{for}\quad 8\leq i \leq 13,\quad \text{and}\quad 16\leq i \leq 21.
\end{equation}

Since $\kappa^L$ is supported only on the identity we must also consider
property~\eqref{PBWconditions-iii} for $g \neq 1$ and find all
$G$-invariant constant $2$-cochains such that $\psi\big(\kappa^C\big)=0$
(i.e., satisfying the mixed Jacobi identity). By~Corollary~\ref{mixedJacobi}, these are
supported on the identity and transpositions and are given by
$\kappa^C=\kappa^C_1+\kappa^C_{\refl}$ with $\kappa^C_1$ defined using $\alpha$
and $\beta$ as in Definition~\ref{def:kappa1} and $\kappa^C_{\refl}$
defined using $c$ as in Definition~\ref{def:kapparefl}.

Now assume $\kappa^C_1+\kappa^C_{\refl}$ clears the first obstruction, i.e.,
that the conditions in~\eqref{ClearFirstObstr} do hold for $\kappa^L_1$, and
consider property~\eqref{PBWconditions-iv} for each of $\kappa^C_1$
and $\kappa^C_{\refl}$. Using the values of $\phi^C_{1,1}(u,v,w)$ and $\phi^C_{g,1}(u,v,w)$ in~Sections~\ref{AllVecInW}--\ref{TwoVecInWTwoIndices} we see that
$\kappa^C_1$ clears the second obstruction for $\kappa^L_1$,
i.e., $\phi\big(\kappa^C_1,\kappa^L_1\big)=0$,
precisely when in addition
\begin{gather}\label{ClearSecondObstr1C}
\gamma_7=\gamma_{15}=0.
\end{gather}
We also see that $\kappa^C_{\refl}$ clears the second obstruction for $\kappa^L_1$,
i.e., that $\phi\big(\kappa^C_{\refl},\kappa^L_1\big)=0$,
precisely when in~addition,
\begin{gather}
\label{ClearSecondObstr1Ref}
\gamma_3=\gamma_6=\gamma_{14}=\gamma_{22}=0.
\end{gather}
Thus $\kappa^L_1$ lifts to $\kappa^L_1+\kappa^C_1+\kappa^C_{\refl}$ if and
only if $\gamma_i=0$ for $1 \leq i \leq 22$.
\end{proof}

\begin{Corollary}\label{SpecialCases}
Theorem~$\ref{thm:LOAMapsDoubledPerm}$ includes precisely the following special cases:
\begin{enumerate}[label=$(\arabic*)$]\itemsep=0pt
\item $\kappa=\kappa^L_1$ satisfying conditions~\eqref{ClearFirstObstr},
\item $\kappa=\kappa^L_1+\kappa^C_1$ satisfying conditions~\eqref{ClearFirstObstr} and~\eqref{ClearSecondObstr1C}, and
\item $\kappa=\kappa^L_1+\kappa^C_1+\kappa^C_{\refl}$ with $\kappa^C_{\refl}\not\equiv 0$
satisfying conditions~\eqref{Obstr2kappaCrefc}--\eqref{Obstr2kappaC1bSimplified}.
\end{enumerate}
\end{Corollary}

\begin{proof}\quad

\textit{Cases~$(1)$ and $(2)$.} These are immediate by the forms of \eqref{ClearSecondObstr1C} and~\eqref{ClearSecondObstr1Ref} since $\kappa^C_1 \equiv 0$ if~and only if $\alpha=\beta=0$ and $\kappa^C_{\refl} \equiv 0$ if and only if $c=0$.

\textit{Case}~(3). Suppose instead that $\kappa^L_1$ lifts to $\kappa^L_1+\kappa^C_1+\kappa^C_{\refl}$
with $\kappa^C_{\refl} \not\equiv 0$. Then
in addition to~\eqref{ClearFirstObstr} and~\eqref{ClearSecondObstr1C},
we have $c \neq 0$. This reduces the conditions in~\eqref{ClearSecondObstr1Ref}
to
\begin{gather}\label{Obstr2kappaCrefc}
a_1=b_6 \phantom{-}\qquad \text{and} \qquad a_1-b_3+b_5=0,\qquad \text{or equivalently}\qquad b_3-b_5-b_6=0,\\
\label{Obstr2kappaCrefd}
b_2=-a_5 \qquad \text{and} \qquad b_2+a_3-a_6=0,\qquad \text{or equivalently}\qquad a_3-a_5-a_6=0.
\end{gather}
These in turn
allow simplification of the conditions for $\kappa^C_1+\kappa^C_{\refl}$ to clear the first obstruction
if~$\kappa^C_{\refl} \not\equiv 0$ clears the second obstruction for $\kappa^L_1$
by reducing~\eqref{ClearFirstObstr} to
\begin{gather}
\label{Obstr1Simplifieda}
a_1(a_4 -a_7) - (b_4-b_7)(a_4+a_6+(n-1)a_7)=0, \\
\label{Obstr1Simplifiedb}
b_1(a_4 - a_7) -(b_4 -b_7)(b_4 + b_6 + (n-1) b_7)=0, \\
\label{Obstr1Simplifiedaprime}
a_2(b_4 -b_7) -(a_4-a_7)(a_4+a_5+(n-1)a_7)=0, \\
\label{Obstr1Simplifiedbprime}
b_2(b_4 -b_7) -(a_4 - a_7)(b_4 + b_5 + (n-1) b_7)=0, \\
\label{Obstr1Simplifiedc}
a_1(a_3 + na_4) + a_5(b_3 + nb_4) -b_1 a_2 - b_5 a_6=0, \\
\label{Obstr1Simplifiedd}
a_1(a_3 + na_7) + a_5(b_3 + nb_7) -b_1 a_2 -b_5 a_6=0, \\
\label{Obstr1Simplifiede}
b_1(a_3 + na_4) + b_5(b_3 + nb_4)- 2b_1 b_2 -2b_5 b_6=0, \\
\label{Obstr1simplifiedf}
b_1(a_3 + na_7) + b_5 (b_3 + nb_7) - 2b_1 b_2 -2b_5 b_6=0, \\
\label{Obstr1Simplifiedg}
-a_2(b_3 + nb_4) + a_6(a_3 + na_4) - 2a_1 a_2 - 2a_5 a_6=0, \\
\label{Obstr1Simplifiedh}
-a_2(b_3 + nb_7) + a_6(a_3 + na_7) - 2a_1 a_2 -2a_5 a_6=0.
\end{gather}
Conditions~\eqref{Obstr2kappaCrefc} and~\eqref{Obstr2kappaCrefd} also reduce the conditions in~\eqref{ClearSecondObstr1C} for $\kappa^C_1$ to
clear the second obstruction for $\kappa^L_1$ assuming $\kappa^C_{\refl} \not\equiv 0$
also clears the second obstruction for $\kappa^L_1$ to
\begin{gather}
\label{Obstr2kappaC1aSimplified}
(\alpha+(n-1)\beta)(b_4-b_7)=0, \\
\label{Obstr2kappaC1bSimplified}
(\alpha+(n-1)\beta)(a_4-a_7)=0.
\end{gather}
Thus
$\kappa^L_1$ lifts to $\kappa^L_1+\kappa^C_1+\kappa^C_{\refl}$ with $\kappa^C_{\refl}\not\equiv 0$
when~\eqref{Obstr2kappaCrefc}--\eqref{Obstr2kappaC1bSimplified} hold.
\end{proof}

\begin{Remark}
Note in Case~(3) of Corollary~\ref{SpecialCases} with $c \neq 0$ that if $\alpha \neq -(n-1)\beta$, then $a_4=a_7$ and $b_4=b_7$.
These combined with $a_3=a_5+a_6$ and $b_3=b_5+b_6$
in~\eqref{Obstr2kappaCrefc}--\eqref{Obstr2kappaCrefd}
mean that the definitions of $\kappa^L_1(x_i,y_j)$ in~\eqref{kappaLxiyj} (but with
$j$ allowed to be $i$) and the definition of~$\kappa^L_1(x_i,y_i)$ in~\eqref{kappaLxiyi}
agree.
Furthermore, conditions~\eqref{Obstr1Simplifieda}--\eqref{Obstr1Simplifiedbprime} then hold and
conditions~\eqref{Obstr1Simplifiedc}--\eqref{Obstr1Simplifiedh},
simplify further to
\begin{gather*}
a_1(a_3 + na_4) + a_5(b_3 + nb_4) - b_1 a_2 - b_5 a_6=0, \\
b_1(a_3 + na_4) + b_5(b_3 + nb_4)- 2b_1 b_2 - 2b_5 b_6=0, \\
-a_2(b_3 + nb_4) + a_6(a_3 + na_4) - 2a_1 a_2 - 2a_5 a_6=0.
\end{gather*}
\end{Remark}

\subsection{Algebraic varieties corresponding to image constraints}\label{algvars}

The homogeneous quadratic PBW conditions $\gamma_i=0$ for $1 \leq i \leq 22$
give rise to a projective variety that controls the parameter space for the
family of maps in Theorem~\ref{thm:LOAMapsDoubledPerm}
and corresponding Lie orbifold algebras in Theorem~\ref{Lieorbifold}.
Based on computations done
for a few specific values of $n$ in Macaulay2~\cite{M2} with the graded reverse
lexicographic monomial ordering and the
parameter order $a_1,\dots,a_7$, $b_1,\dots,b_7$, $\alpha$, $\beta$, $c$
we conjecture that this projective variety has dimen\-sion seven. It contains the
lattice of subvarieties described in Table~\ref{table:subvars} arising from
the PBW conditions and additional constraints on how $\im \kappa^L_1$ relates to
subrepresentations of the doubled permutation representation.
While not needed for Theorem~\ref{thm:LOAMapsDoubledPerm}, these are included as being of potential independent interest.

\begin{Remark}\label{rmk:subrep1}
Observe that
$x_i=\bar{x}_i+\tfrac{1}{n}x_{[n]}$ and $y_i=\bar{y}_i+\tfrac{1}{n}y_{[n]}$
can be used in~\eqref{kappaLxx}--\eqref{kappaLxiyj} to decompose the results
according to $V\cong\mathfrak{h}^* \oplus \iota^* \oplus \mathfrak{h} \oplus \iota$:
\begin{gather}
\kappa^L_1(x_i,x_j)=a_1(\bar{x}_i-\bar{x}_j)+b_1(\bar{y}_i-\bar{y}_j), \\
\kappa^L_1(y_i,y_j)=a_2(\bar{x}_i-\bar{x}_j)+b_2(\bar{y}_i-\bar{y}_j), \\
\kappa^L_1(x_i,y_i)=a_3\bar{x}_i+\frac{1}{n}(a_3+na_4)\n{x}+b_3 \bar{y}_i+\frac{1}{n}(b_3+nb_4)\n{y}, \\
\kappa^L_1(x_i,y_j)=a_5 \bar{x}_i\!+a_6 \bar{x}_j\!+\frac{1}{n}(a_5\!+a_6\!+na_7)\n{x}\!+b_5 \bar{y}_i\!+b_6 \bar{y}_j\!+\frac{1}{n}(b_5\!+b_6\!+nb_7)\n{y}.
\end{gather}
\end{Remark}

Using Remark~\ref{rmk:subrep1} we see that if the image of $\kappa^L_1$ is contained within $\mathfrak{h}^*\oplus \mathfrak{h}$
then the coefficients of $x_{[n]}$ and $y_{[n]}$ are zero;
i.e.,
\begin{alignat}{3}
\label{DoubledStdRepa}
& a_3+na_4=0, \qquad && b_3+nb_4=0,& \\
\label{DoubledStdRepb}
& a_5+a_6+na_7=0, \qquad && b_5+b_6+nb_7=0.&
\end{alignat}
If the image of $\kappa^L_1$ contains at most one copy of the trivial representation,
i.e., some linear combination of $x_{[n]}$ and $y_{[n]}$, then the coefficients of
$x_{[n]}$ and $y_{[n]}$ satisfy this weaker condition:
\begin{gather}\label{TrivialRep}
(a_3+na_4)(b_5+b_6+nb_7)=(b_3+nb_4)(a_5+a_6+na_7).
\end{gather}
Similarly, if the image of $\kappa^L_1$ contains at most one copy of the standard representation
then the coefficients of $\bar{x_i}$ and $\bar{y_i}$ satisfy these conditions:
\begin{gather}\label{StdRep}
a_i b_j=b_i a_j \qquad \text{for}\quad 1\leq i<j \leq 6\quad \text{and}\quad i,j \neq 4.
\end{gather}
Lastly, if the image of $\kappa^L_1$ is contained within $\iota^* \oplus \iota$
then the coefficients of $\bar{x_i}$ and $\bar{y_i}$ are zero:
\begin{gather}\label{DoubledTrivialRep}
a_i=b_i=0 \qquad \text{for}\quad 1\leq i \leq 6 \quad\text{and}\quad i \neq 4.
\end{gather}

In Table~\ref{table:subvars} we list constraints on $\im \kappa^L_1$,
resulting conditions, in addition to \eqref{ClearFirstObstr}--\eqref{ClearSecondObstr1Ref},
needed to describe the subvariety of maps subject to each constraint,
and the conjectured dimension and degree based on computations.
Dropping condition~\eqref{TrivialRep} in the first three rows, i.e., allowing an additional summand of the trivial representation in $\im \kappa^L_1$ in those cases, did not change the computed dimension or degree of the variety.

\begin{table}[hbt]
\begin{center}
\caption{Conjectured dimension and degree of projective varieties 
for Lie orbifold algebras obtained from constraints on $\im \kappa^L_1$. When $\im \kappa^L_1$ is contained in the trivial representation, the resulting seven defining polynomials, the dimension and degree, and related simple $\kappa$ maps can be found by hand.}\label{table:subvars}

\vspace{2mm}

\renewcommand{\arraystretch}{1.2}%
\begin{tabular}{l|l|c|c}
\hline
\multicolumn{1}{c|}{Constraint on $\im \kappa^L_1$} & \multicolumn{1}{c|}{Conditions} & \multicolumn{1}{c|}{Dimension} & \multicolumn{1}{c}{Degree}\\
\hline
contained in $\mathrm{perm} \oplus \mathrm{std}$ & \eqref{TrivialRep} & 7 & 8 \\
contained in $\mathrm{perm}$ & \eqref{TrivialRep}, \eqref{StdRep} & 6 & 30 \\
contained in $\mathrm{triv}$ & \eqref{TrivialRep}, \eqref{StdRep}, \eqref{DoubledTrivialRep} & 4 & 1 \\
contained in $\mathrm{std} \oplus \mathrm{std}$ & \eqref{DoubledStdRepa}, \eqref{DoubledStdRepb}, \eqref{TrivialRep} & 6 & 6 \\
contained in std & \eqref{DoubledStdRepa}, \eqref{DoubledStdRepb}, \eqref{TrivialRep}, \eqref{StdRep} & 5 & 4
\\
\hline
\end{tabular}
\end{center}
\end{table}

\section[Drinfeld orbifold algebra maps that deform S(W* + W) Sn]
{Drinfeld orbifold algebra maps that deform $\boldsymbol{S(W^* \oplus W)\# S_n}$}\label{Lifting}

In Section~\ref{InvSkewSymmBilMaps} we defined a pre-Drinfeld orbifold algebra map
$\kappa^L_{\refl}$. Here we aim to lift $\kappa^L_{\refl}$ to a Drinfeld orbifold
algebra map. We~evaluate $\phi\big(\kappa^L_{\refl},\kappa^L_{\refl}\big)$
and clear the first obstruction by defining a $G$-invariant $2$-cochain
$\kappa^C_{\tri}$ such that $\phi\big(\kappa^L_{\refl},\kappa^L_{\refl}\big)=2\psi\big(\kappa^C_{\tri}\big)$. We~then clear the second obstruction by showing
$\phi\big(\kappa^C_{\refl}+\kappa^C_{\tri},\kappa^L_{\refl}\big)=0$ and giving conditions on
the parameters~$\alpha$ and~$\beta$ for $\kappa^C_1$ and $a$, $a^\perp$, $b$, $b^\perp$
for $\kappa^L_{\refl}$ such that $\phi\big(\kappa^C_1,\kappa^L_{\refl}\big)=0$. It follows
in Theorem~\ref{thm:DOAMapsDoubledPerm} that
$\kappa^L_{\refl}+\kappa^C_{\refl}+\kappa^C_{\tri}$ is always a Drinfeld orbifold
algebra map and $\kappa^L_{\refl}+\kappa^C_1+\kappa^C_{\refl}+\kappa^C_{\tri}$ is
a~Drinfeld orbifold algebra map precisely when conditions~\eqref{Obstr2PhiC1C2C3L2=0a}
and~\eqref{Obstr2PhiC1C2C3L2=0b} hold.

Characterizing in general when $\kappa^L_1+\kappa^L_{\refl}$ lifts is
straightforward but rather more involved. Instead, in
Theorem~\ref{thm:DOAMapsCombined} we provide a nontrivial choice of parameters
and verify in that case that it is possible to lift simultaneously to
$\kappa^L_1+\kappa^L_{\refl}+\kappa^C_1+\kappa^C_{\refl}+\kappa^C_{\tri}$.
There could very well be other parameter choices for successfully lifting $\kappa_1^L+\kappa^L_{\refl}$. This is indicated by the question marks in Table~\ref{SummaryDOAMapsDoubledPerm} summarizing the results in this
and the previous section.

\subsection[Clearing obstructions to deformations of S(W* + W) Sn]
{Clearing obstructions to deformations of $\boldsymbol{S(W^* \oplus W)\# S_n}$}\label{ObstructionDetails}

We begin by recalling a lemma from~\cite{FGK} that allows
for a reduction in the computations necessary to remove obstructions and lift $\kappa$ maps.

\medskip\noindent
{\bf Invariance relations.}
Recall that a cochain $\mu=\!\sum_{g\in G}\mu_gg$ with components $\mu_g\colon \bigwedge^k V\!\to S(V)$ is $G$-invariant if and only if $\actby{h}\mu_g=\mu_{hgh^{-1}}$ for all $g,h\in G$.
Equivalently,
\begin{gather*}
\actby{h}(\mu_g(v_1,\ldots,v_k))=\mu_{hgh^{-1}}\big(\actby{h}v_1,\ldots,\actby{h}v_k\big)
\end{gather*}
for all $g,h\in G$ and $v_1,\ldots,v_k\in V$. Thus a $G$-invariant cochain is
determined by its components for a set of conjugacy class representatives.

In the following lemma that applies to any finite group, one can let $\mu=\kappa^L$ or $\mu=\kappa^C$ and let $\nu=\kappa^L$ to see that if $\kappa^L$ and $\kappa^C$ are $G$-invariant,
then $\phi\big(\kappa^{*},\kappa^L\big)$ and $\psi(\kappa^{*})$ are also $G$-invariant. This is helpful because, for instance, if $\phi_g=2\psi_g$ for some $g\in G$, then acting by~$h\in G$ on~both sides shows $\phi_{hgh^{-1}}=2\psi_{hgh^{-1}}$ also. Thus if $\phi_g=2\psi_g$ for all $g$ in a set of conjugacy class representatives, then $\phi\big(\kappa^L,\kappa^L\big)=2\psi\big(\kappa^C\big)$. Similar reasoning applies to properties~\eqref{PBWconditions-ii} and~\eqref{PBWconditions-iv} of a Drinfeld orbifold algebra map.

\begin{Lemma}[{\rm \cite[Lemma 5.1]{FGK}}]	\label{le:orbitphixyphig}
	Let $G$ be a finite group acting linearly on $V\cong\mathbb{C}^n$. If~$\mu$ and~$\nu$ are $G$-invariant $2$-cochains with $\nu$ linear and $\mu$ linear or constant,
	then $\phi(\mu,\nu)$ and $\psi(\mu)$ are $G$-invariant. Specifically, at
	the component level, for all $x,y,h\in G$ and $v_1,v_2,v_3\in V$ we have
	\begin{gather*}
	\actby{h}(\phi_{x,y}(v_1,v_2,v_3))=\phi_{hxh^{-1},hyh^{-1}}\big(\actby{h}v_1,\actby{h}v_2,\actby{h}v_3\big).
	\end{gather*}
\end{Lemma}

\subsection{Clearing the first obstruction}

We begin by recording simplifications of a summand, $\phi^{*}_{\sigma,\tau}$,
of the component $\phi_{g}^{*}$ of $\phi\big(\kappa^{*}_{\refl},\kappa^L_{\refl}\big)$,
where $*$ stands for $L$ or $C$.
Simplification of $\phi^{*}_{\sigma,\tau}(u,v,w)$ depends on
the location of the basis vectors relative to $W$ and $W^*$ and relative to
the fixed spaces $V^{\sigma}$ and $V^{\tau}$, so recall that $v^*$ is the
vector dual to $v$ and define the following indicator function. For~$g \in S_n$ and $v \in V$, let
\begin{gather*}
\delta_g(v)=\begin{cases}
1 & \text{if}\quad v \in V^{g}, \\
0 & \text{otherwise}.
\end{cases}
\end{gather*}
Note that for $g \in G$ and $v \in V$,
\begin{gather*}
\delta_g(v^*)=\delta_g(v).
\end{gather*}

\noindent\textbf{Remark.}
Let $\phi^*_{\sigma,\tau}$ be as in Lemma~\ref{lemma:simplification*2L2}.
Then for all $u,v,w \in V$ we have
\begin{equation}
\label{tauinv}
\phi^*_{\sigma,\tau}(\actbylessspace{\tau}u,\actbylessspace{\tau}v,\actbylessspace{\tau}w)
=\phi^*_{\sigma,\tau}(u,v,w).
\end{equation}
This follows from the definition of $\phi^*_{\sigma,\tau}$ and that $\kappa^L_{\tau}$
is $\tau$-invariant.

\begin{Lemma}\label{lemma:simplification*2L2}
Let $\kappa^*_{\refl}$ with $*=L$ or $*=C$ be as in Definition~$\ref{def:kapparefl}$,
with common parameters $a, a^\perp, b, b^\perp \in \mathbb{C}$.
Denote a term of the component $\phi^*_g$ of $\phi\big(\kappa^*_{\refl},\kappa^L_{\refl}\big)$
by $\phi^*_{\sigma,\tau}$, where $\sigma$ and $\tau$ are transpositions
such that $\sigma \tau = g$.
	\begin{enumerate}[label={$\arabic*.$}]\itemsep=0pt
		\item[$(1)$] 
If $u,v,w \in W$, $u,v,w \in W^*$, $u,v\in V^{\tau}$, or
		$u \in V^{\tau}\cap V^{\sigma}$, then $\phi^{*}_{\sigma,\tau}(u,v,w)=0$.
		\item[$(2)$] 
If $u \in V^{\tau} \setminus V^{\sigma}$ and $v \notin V^\tau$, then the
		basis vectors moved by $\tau$ are of the form $v$, $\actby{\tau}v$, $v^*$,~$\actby{\tau}v^*$.
		We have $\phi^*_{\sigma,\tau}(v,\actbylessspace{\tau}v,u^*)=0$,
		\begin{gather*}
		\phi^*_{\sigma,\tau}(u,v,v^*)
		=-\phi^*_{\sigma,\tau}(u,v,\actbylessspace{\tau}v^*)
		=\begin{cases}
		2\big(a^\perp-a\big)[1-\delta_{\tau}(\actbylessspace{\sigma}u)]\kappa^*_{\sigma}(u^*,u)
& \text{if}\quad u \in W, \\
		2\big(b^\perp-b\big)[1-\delta_{\tau}(\actbylessspace{\sigma}u)]\kappa^*_{\sigma}(u,u^*)
& \text{if}\quad u \in W^*,
		\end{cases}
		\end{gather*}
		and by~\eqref{tauinv},
	\begin{gather*}
\phi^*_{\sigma,\tau}(u,\actby{\tau}v,\actby{\tau}v^*)=-\phi^*_{\sigma,\tau}(u,\actby{\tau}v,v^*) =\phi^*_{\sigma,\tau}(u,v,v^*).
\end{gather*}
		\item[$(3)$] 
If $v \notin V^{\tau}$, then
		\begin{gather*}
		\phi^*_{\sigma,\tau}(v,\actby{\tau}v,v^*)
		=\begin{cases}
		2\big(a^\perp-a\big)\big[\delta_{\tau}(\actbylessspace{\sigma}v)\kappa^*_{\sigma}(v^*,v) +\delta_{\tau}(\actby{\sigma\tau}v)\kappa^*_{\sigma}(\actbylessspace{\tau}v^*,\actbylessspace{\tau}v) \big] & \text{if}\quad v \in W,
\\
2\big(b^\perp-b\big)\big[\delta_{\tau}(\actbylessspace{\sigma}v)\kappa^*_{\sigma}(v,v^*) +\delta_{\tau}(\actby{\sigma\tau}v)\kappa^*_{\sigma}(\actbylessspace{\tau}v,\actbylessspace{\tau}v^*) \big] & \text{if}\quad v \in W^*,
		\end{cases}
		\end{gather*}
		and by~\eqref{tauinv},
		\begin{gather*}
		\phi^*_{\sigma,\tau}(v,\actbylessspace{\tau}v,\actbylessspace{\tau}v^*)
		=-\phi^*_{\sigma,\tau}(v,\actbylessspace{\tau}v,v^*).
		\end{gather*}
	\end{enumerate}
\end{Lemma}

\begin{proof} \quad

\textit{Case}~(1). When $u,v,w \in W$ or $u,v,w \in W^*$, then $\phi^{*}_{\sigma,\tau}(u,v,w)=0$ because
	\begin{gather}
	\label{phist}
	\phi^{*}_{\sigma,\tau}(u,v,w)=\phibulletxy{\sigma}{\tau}{u}{v}{w},
	\end{gather}
and $\kappa^L_{\tau}$ is zero whenever both
input vectors are in $W^*$ or both are in $W$. As in~\cite{FGK},
$\phi^{*}_{\sigma,\tau}(u,v,w)=0$ when $u, v \in V^{\tau}$ follows from~\eqref{phist}
and that $V^{\tau}\subseteq \ker\kappa^*_{\tau}$, while $\phi^{*}_{\sigma,\tau}(u,v,w)=0$
when $u \in V^{\tau} \cap V^{\sigma}$ uses also $V^{\sigma}\subseteq \ker\kappa^*_{\sigma}$.

\textit{Case}~(2). Assume $u \in V^{\tau} \setminus V^{\sigma}$ and $v \notin V^\tau$.
First note that
	\begin{gather*}
	\phi^*_{\sigma,\tau}(v,\actbylessspace{\tau}v,u^*)=0
	\end{gather*}
by using~\eqref{tauinv}
and the alternating property to see that $\phi^*_{\sigma,\tau}(v,\actbylessspace{\tau}v,u^*)=-\phi^*_{\sigma,\tau}(v,\actbylessspace{\tau}v,u^*)$.
		
We can reduce to
$\phi^*_{\sigma,\tau}(u,v,v^*)=2\kappa^*_{\sigma}\big(u,\kappa^L_{\tau}(v,v^*)\big)$
by using $u \in V^{\tau} \subseteq \ker\kappa^*_{\tau}$ in~\eqref{phist}.
Using bilinearity and $V^{\sigma} \subseteq \ker \kappa^*_{\sigma}$, the right hand
side is a linear combination of expressions $\kappa^*_{\sigma}\big(u,\actby{h}u\big)$ and
$\kappa^*_{\sigma}\big(u,\actby{h}u^*\big)$ for $h \in \langle \sigma \rangle$.
The appropriate coefficients in terms of $a$, $a^\perp$, $b$, $b^\perp$ can be described
in terms of the indicator function for the fixed space of $\tau$ and depend on
whether $u \in W^*$ or~$u \in W$.
Also note that $\sum_{h \in \langle \sigma \rangle} \actby{h}u \in V^{\sigma} \subseteq \ker \kappa^*_{\sigma}$.
Thus, for $u \in W^*$ we have
\begin{gather*}
\kappa^*_{\sigma}\big(u,\kappa^L_{\tau}(v,v^*)\big)
\\ \quad
{}=\sum\limits_{h \in \langle \sigma \rangle} \!\! \big[a\big(1-\delta_{\tau}\big(\actby{h}u\big)\big)+a^\perp\delta_{\tau}\big(\actby{h}u\big)\big]\kappa^*_{\sigma}\big(u,\actby{h}u\big) +\big[b\big(1-\delta_{\tau}\big(\actby{h}u^*\big)\big)+b^\perp\delta_{\tau}\big(\actby{h}u^*\big)\big] \kappa^*_{\sigma}(u,\actby{h}u^*)
\\ \quad
{}=\big(a^\perp-a\big)\sum\limits_{h \in \langle \sigma \rangle}
\delta_{\tau}\big(\actby{h}u\big)\kappa^*_{\sigma}\big(u,\actby{h}u\big)
+\big(b^\perp-b\big)\sum\limits_{h \in \langle \sigma \rangle}
\delta_{\tau}\big(\actby{h}u^*\big)\kappa^*_{\sigma}\big(u,\actby{h}u^*\big)
\\ \quad
{}=\big(b^\perp-b\big)\big[\delta_{\tau}(u)\kappa^*_{\sigma}(u,u^*)
+\delta_{\tau}(\actbylessspace{\sigma}u)\kappa^*_{\sigma}(u,\actbylessspace{\sigma}u^*)\big],
\end{gather*}
where the term with coefficient $a^\perp-a$ is zero because $\kappa^*_{\sigma}(u,u)=\kappa^*_{\sigma}(u,\actbylessspace{\sigma}u)=0$.
Since $u \in V^{\tau}$ and $\kappa^*_{\sigma}(u,\actbylessspace{\sigma}u^*)=-\kappa^*_{\sigma}(u,u^*)$, this yields
\begin{gather*} \phi^*_{\sigma,\tau}(u,v,v^*)=2\kappa^*_{\sigma}\big(u,\kappa^L_{\tau}(v,v^*)\big)=2\big(b^\perp-b\big) \big[1-\delta_{\tau}(\actbylessspace{\sigma}u)\big]\kappa^*_{\sigma}(u,u^*).
\end{gather*}
	
When $u \in W$, the calculation of $\kappa^*_{\sigma}\big(u,\kappa^L_{\tau}(v,v^*)\big)$
involves a sign difference, and the coefficients on the two sums
are reversed, so the one with coefficient $a^\perp-a$ survives and yields
$\phi^*_{\sigma,\tau}(u,v,v^*)=2\big(a^\perp-a\big)\big[1-\delta_{\tau}(\actbylessspace{\sigma}u)\big]\kappa^*_{\sigma}(u^*,u)$.

Similar calculations show
\begin{gather*}
\phi^*_{\sigma,\tau}(u,v,\actbylessspace{\tau}v^*)=-\phi^*_{\sigma,\tau}(u,v,v^*).
\end{gather*}
		
\textit{Case}~(3). Assume $v \notin V^{\tau}$.
Then $\kappa^L_{\tau}(v,\actbylessspace{\tau}v)=0$ and hence
\begin{gather*}
\phi^*_{\sigma,\tau}(v,\actbylessspace{\tau}v,v^*) =\kappa^*_{\sigma}\big(v+\actbylessspace{\tau}v,\kappa^L_{\tau} (\actbylessspace{\tau}v,v^*)\big)+\kappa^*_{\sigma}\big(\actbylessspace{\tau}v+v,\kappa^L_{\tau}(v^*,v)\big)
=-2\kappa^*_{\sigma}\big(v+\actbylessspace{\tau}v,\kappa^L_{\tau}(v,v^*)\big).
\end{gather*}
A calculation as in case~(2), using
$\kappa^L_\tau(\actbylessspace{\tau}v,\actbylessspace{\tau}v^*)=\kappa^L_{\tau}(v,v^*)$
and that $\delta_{\tau}(v)=\delta_{\tau}(\actbylessspace{\tau}v)=0$ yields
\begin{gather*}
\phi^*_{\sigma,\tau}(v,\actbylessspace{\tau}v,v^*)
=-2\left[\kappa^*_{\sigma}\big(v,\kappa^L_{\tau}(v,v^*)\big) +\kappa^*_{\sigma}\big(\actbylessspace{\tau}v,\kappa^L_{\tau} (\actbylessspace{\tau}v,\actbylessspace{\tau}v^*)\big)\right]
\\ \hphantom{\phi^*_{\sigma,\tau}(v,\actbylessspace{\tau}v,v^*) }
=\begin{cases}
2\big(a^\perp-a\big)\big[\delta_{\tau}(\actbylessspace{\sigma}v)\kappa^*_{\sigma}(v^*,v) +\delta_{\tau}(\actby{\sigma\tau}v)\kappa^*_{\sigma}(\actbylessspace{\tau}v^*,\actbylessspace{\tau}v) \big] & \text{if}\quad v \in W,
\\
2\big(b^\perp-b\big)\big[\delta_{\tau}(\actbylessspace{\sigma}v)\kappa^*_{\sigma}(v,v^*) +\delta_{\tau}(\actby{\sigma\tau}v)\kappa^*_{\sigma}(\actbylessspace{\tau}v,\actbylessspace{\tau}v^*) \big] & \text{if}\quad v \in W^*.
\end{cases}	
\tag*{\qed}
\end{gather*}
\renewcommand{\qed}{}
\end{proof}

As mentioned in the outline of the proof of Theorem~\ref{thm:DOAMapsDoubledPerm}, the next two propositions are used to evaluate both $\phi\big(\kappa^L_{\refl},\kappa^L_{\refl}\big)$ and $\phi\big(\kappa^C_{\refl},\kappa^L_{\refl}\big)$.

\begin{Proposition}\label{prop:zerocases}
	Let $\kappa^{*}_{\refl}$ with $*=L$ or $*=C$ be as in Definition~$\ref{def:kapparefl}$.
	For $g \in S_n$ where $n \geq 3$, let $\phi^{*}_g$ be the $g$-component of $\phi\big(\kappa^{*}_{\refl},\kappa^L_{\refl}\big)$. If~$g$ is not a $3$-cycle then $\phi^{*}_g \equiv 0$.
\end{Proposition}

\begin{proof}
Since $\kappa^{*}_{\refl}$ is supported only on transpositions and the only cycle
types that arise as a~product of two transpositions are the identity, double transpositions,
and $3$-cycles, it suffices to consider only the components $\phi_1^*$ and $\phi_g^*$
where $g$ is a double transposition, and in fact only~$\phi^{*}_1$ when $n=3$.
Since
$\actby{h}\sigma\actby{h}\tau=\actby{h}(\sigma \tau)$, it suffices to use only
representatives of orbits of factor pairs under the action of $S_n$ by diagonal conjugation.
	
\textit{Case}~1 $(g=1)$.
The identity component $\phi^{*}_1$ of $\phi(\kappa^{*}_{\refl},\kappa^L_{\refl})$ is
a sum of terms $\phi^{*}_{\sigma,\sigma^{-1}}$, where $\sigma$ ranges over the set of transpositions in~$S_n$. For~each transposition $\sigma$, since $\sigma^{-1}=\sigma$ we have
$\im\kappa^L_{\sigma^{-1}}=\im\kappa^L_{\sigma}\subseteq V^{\sigma} \subseteq\ker\kappa^{*}_{\sigma}$,
and thus $\kappa^{*}_{\sigma}(u,\kappa^L_{\sigma^{-1}}(v,w))=0$ for all $u,v,w\in V$.
It follows that $\phi^{*}_{\sigma,\sigma^{-1}}\equiv0$ for each transposition $\sigma\in S_n$,
and hence, $\phi^{*}_1\equiv 0$.

\textit{Case}~2 $(g=(12)(34))$.
Note that $\phi^{*}_g=
\phi^{*}_{(12),(34)}+\phi^{*}_{(34),(12)}$.
Since we know that both $\im\kappa^L_{(34)}\subseteq V^{(12)} \subseteq\ker\kappa^{*}_{(12)}$ and
$\im\kappa^L_{(12)}\subseteq V^{(34)} \subseteq\ker\kappa^{*}_{(34)}$, we see
that $\phi^{*}_g(u,v,w)\equiv 0$.
\end{proof}

\begin{Proposition}\label{prop:3-cycle}
Let $\kappa_{\refl}=\kappa^L_{\refl}+\kappa^C_{\refl}$ be as in
Definition~$\ref{def:kapparefl}$ with parameters $a,b,c\in\mathbb{C}$, and let
$\phi_g^{*}$ denote the $g$-component
of $\phi\big(\kappa^*_{\refl},\kappa^L_{\refl}\big)$, where $*=L$
or $*=C$ and $g$ is a $3$-cycle. Then $\phi_g^C\equiv 0$. For~$\phi_g^L$ we have
$\phi_g^L(u,v,w)=0$ if $u\in V^g$,
and for $v \notin V^g$, $\phi_g^L\big(v,\actby{g}v,\actbylessspace{g^2}v\big)=0$, and
\begin{align*} \phi_{g}^L(v,\actby{g}v,v^*)&=\phi_{g}^L\big(\actby{g}v,\actbylessspace{g^2}v,v^*\big)=\phi_{g}^L\big(\actbylessspace{g^2}v,v,v^*\big) \\
&=\begin{cases}
2\big(a^\perp-a\big)\big(b^\perp-b\big)(\actby{g}v-v)+2\big(a^\perp-a\big)^2(\actby{g}v^*-v^*) & \text{if\quad $v \in W$,} \\
2\big(a^\perp-a\big)\big(b^\perp-b\big)(\actby{g}v-v)+2\big(b^\perp-b\big)^2(\actby{g}v^*-v^*) & \text{if\quad $v \in W^*$,}
\end{cases}
\end{align*}
with values on triples involving a third basis vector of the form
$\actby{g}v^*$ or $\actbylessspace{g^2}v^*$ obtained by acting by~$g$ or $g^2$ respectively.
\end{Proposition}

\begin{proof}
By the orbit property in Lemma~\ref{le:orbitphixyphig}, it suffices to evaluate $\phi^*_g$ for the conjugacy
class representative $g=(123)$. Note that $Z(g)=\langle(123)\rangle\times \Sym_{\{4,\dots,n\}}$,
and the factorizations of $g$ as a product of transpositions are all in the same $Z(g)$-orbit
under diagonal conjugation, so
\begin{gather*}
\phi_g^{*}=\phi_{(12),(23)}^{*}+\phi_{(23),(31)}^{*}+\phi_{(31),(12)}^{*}.
\end{gather*}
For each pair of transpositions $\sigma$ and $\tau$ with $\sigma\tau=g=(123)$,
we have $V^g \subseteq V^\sigma \cap V^\tau$ and
$V^\sigma \cap V^\tau \subseteq \ker \phi^*_{\sigma,\tau}$
by Lemma~\ref{lemma:simplification*2L2}(1). So if any vector in a basis triple is in $V^g$, then $\phi^*_g$
evaluates to zero. It remains only to consider triples
$\{u,v,w\} \subseteq \{x_1,x_2,x_3,y_1,y_2,y_3\}$.

Lemma~\ref{lemma:simplification*2L2}(1) yields
$\phi_g^{*}(x_1,x_2,x_3)=\phi_g^{*}(y_1,y_2,y_3)=0$.
There are, up to permutation, $\frac{1}{2}\binom{6}{3}-1=9$
basis triples with two elements in $W^*$ and one element in $W$:
\begin{alignat*}{3}
&x_1, x_2, y_1, &\qquad
&x_2, x_3, y_1, &\qquad
&x_3, x_1, y_1, \\
&x_2, x_3, y_2, &\qquad
&x_3, x_1, y_2, &\qquad
&x_1, x_2, y_2, \\
&x_3, x_1, y_3, &\qquad
&x_1, x_2, y_3, &\qquad
&x_2, x_3, y_3.
\end{alignat*}
The $G$-invariance in Lemma~\ref{le:orbitphixyphig} yields
	\begin{gather*}
	\phi^*_g(\actby{g}u,\actby{g}v,\actby{g}w)=\actby{g}\phi^*_g(u,v,w),
	\end{gather*}
and thus since each of the three columns of basis triples is a $g$-orbit,
it suffices to compute $\phi^*_g$ on just the basis triples
	\begin{gather*}
	x_1,x_2,y_1, \qquad x_2,x_3,y_1, \qquad \text{ and} \qquad x_3,x_1,y_1.
	\end{gather*}
For each of these three basis triples, $u$, $v$, $w$, the result $\phi^*_g(u,v,w)$
will be the same by Lemma~\ref{lemma:simplification*2L2} (although the reason
varies for a given term on different triples), namely
\begin{align*}
\phi^*_g(u,v,w)
&=\phi_{(12),(23)}^{*}(u,v,w)
+\phi_{(23),(31)}^{*}(u,v,w)
+\phi_{(31),(12)}^{*}(u,v,w) \\
&=0-2\big(b^\perp-b\big)\kappa^{*}_{(23)}(x_3,y_3)+2\big(b^\perp-b\big)\kappa^{*}_{(31)}(x_1,y_1) \\
&=
\begin{cases}
0 & \text{if}\quad $*=C$, \\
-2\big(b^\perp-b\big)\big[ax_{2,3} + a^\perp x_{2,3}^\perp + by_{2,3} + b^\perp y_{2,3}^\perp\big] \\ \qquad
{}+2\big(b^\perp-b\big)\big[ax_{1,3} + a^\perp x_{1,3}^\perp + by_{1,3} + b^\perp y_{1,3}^\perp\big] &
\text{if}\quad *=L
\end{cases} \\
&=\begin{cases}
0 & \text{if}\quad *=C, \\
2\big(b^\perp-b\big)\big[\big(a^\perp-a\big)(x_2-x_1) + \big(b^\perp-b\big)(y_2-y_1)\big] & \text{if}\quad *=L.
\end{cases}
\end{align*}
A similar reduction and computation applies to the nine basis triples with two elements
in $W$ and one element in $W^*$, and that combined with the orbiting properties in
Lemma~\ref{le:orbitphixyphig} lead to the conclusion in the statement.
\end{proof}

Using the form of $\phi\big(\kappa^L_{\refl},\kappa^L_{\refl}\big)$ in Propositions~\ref{prop:zerocases} and~\ref{prop:3-cycle} we define
the cochain $\kappa^C_{\tri}$ in~Defi\-nition~\ref{def:kappatri} to
ensure $\phi\big(\kappa^L_{\refl},\kappa^L_{\refl}\big)=2\psi\big(\kappa^C_{\tri}\big)$
as in the next proposition.

\begin{Proposition}\label{prop:phiL2L2=2psiC3}
 Let $\kappa^L_{\refl}$ and $\kappa^C_{\tri}$
 be as in Definitions~$\ref{def:kapparefl}$ and $\ref{def:kappatri}$, with common
 parameters $a,a^\perp,b,b^\perp \in\mathbb{C}$. Then
 $\phi\big(\kappa^L_{\refl},\kappa^L_{\refl}\big)=2\psi\big(\kappa^C_{\tri}\big)$.
\end{Proposition}

\begin{proof}
	We compare the component $\phi_g$ of $\phi\big(\kappa^L_{\refl},\kappa^L_{\refl}\big)$
	with the component $2\psi_g$ of $2\psi\big(\kappa^C_{\tri}\big)$. If~$g$ is not a $3$-cycle, then $\phi_g\equiv 0$ by Proposition~\ref{prop:zerocases};
	and $\kappa^C_{\tri}$ is not supported on $g$, so $\psi_g\equiv 0$ as well.
	If~$g$ is a $3$-cycle, then it suffices to compare components $\phi_g$ and $2\psi_g$
	on basis triples of the form in the statement of Proposition~\ref{prop:3-cycle}.
	
\textit{Case}~1. If $u\in V^g$, then
$\phi_g(u,v,w)=0$
by Proposition~\ref{prop:3-cycle} and $\psi_g(u,v,w)=0$ because $\actby{g}u-u=0$ and $V^g\subseteq\ker\kappa^C_g$.

\textit{Case}~2. If $v \notin V^g$, then $\phi_{g}\big(v,\actby{g}v,\actbylessspace{g^2}v\big)=0$
by Proposition~\ref{prop:3-cycle} and $\psi_{g}\big(v,\actby{g}v,\actbylessspace{g^2}v\big)=0$
by $g$-invariance of $\kappa^C_g$.
	
\textit{Case}~3. For triples of the form
$(v,\actby{g}v,v^*)$, $\big(\actby{g}v,\actbylessspace{g^2}v,v^*\big)$, and $\big(\actbylessspace{g^2}v,v,v^*\big)$ with $v \in V^g$, use Proposition~\ref{prop:3-cycle} to find $\phi_{g}(v,\actby{g}v,v^*)=\phi_{g}\big(\actby{g}v,\actbylessspace{g^2}v,v^*\big)=\phi_{g}\big(\actbylessspace{g^2}v,v,v^*\big)$
and Definition~\ref{def:kappatri} to confirm that $\phi_g=2\psi_g$ on each such triple, using that
\begin{gather*}
	\psi_{g}(v,\actby{g}v,v^*)
	=\kappa^C_g(v,\actby{g}v)(\actby{g}v^*-v^*)
	+\kappa^C_g(\actby{g}v,v^*)(\actby{g}v-v)
	+\kappa^C_g(v^*,v)\big(\actbylessspace{g^2}v-\actby{g}v\big), \\
	\psi_{g}\big(\actby{g}v,\actbylessspace{g^2}v,v^*\big)
	=\kappa^C_g\big(\actby{g}v,\actbylessspace{g^2}v\big)(\actby{g}v^*-v^*)
	+\kappa^C_g\big(\actbylessspace{g^2}v,v^*\big)\big(\actbylessspace{g^2}v-\actby{g}v\big)
	+\kappa^C_g(v^*,\actby{g}v)\big(v-\actbylessspace{g^2}v\big), \\
	\psi_{g}\big(\actbylessspace{g^2}v,v,v^*\big)
	=\kappa^C_g\big(\actbylessspace{g^2}v,v\big)(\actby{g}v^*-v^*)
	+\kappa^C_g(v,v^*)\big(v-\actbylessspace{g^2}v\big)
	+\kappa^C_g\big(v^*,\actbylessspace{g^2}v)(\actby{g}v-v\big).
\tag*{\qed}
\end{gather*}
\renewcommand{\qed}{}
\end{proof}

\subsection{Clearing the second obstruction}

The final step in determining when the cochain
$\kappa=\kappa^L_{\refl}+\kappa^C_{\tri}+\kappa^C_{\refl}+\kappa^C_1$
is a Drinfeld orbifold algebra map is to understand when
$\phi\big(\kappa^C_{\tri}+\kappa^C_{\refl}+\kappa^C_1,\kappa^L_{\refl}\big)=0$,
which is stated as Corollary~\ref{cor:phiC3L2=0} and follows immediately from
Propositions~\ref{prop:zerocases} and~\ref{prop:3-cycle} and Lemmas~\ref{lemma:simplificationC3L2}
and~\ref{lemma:simplification*1L2}. This clears the second obstruction and
completes the proof of Theorem~\ref{thm:DOAMapsDoubledPerm}.

\begin{Lemma}\label{lemma:simplificationC3L2}
Let $\kappa^L_{\refl}$ and $\kappa^C_{\tri}$ be as in Definitions~$\ref{def:kapparefl}$ and~$\ref{def:kappatri}$,
with common parameters $a, a^\perp, b, b^\perp \in \mathbb{C}$.
Denote a term of the component $\phi_g$ of $\phi\big(\kappa^C_{\tri},\kappa^L_{\refl}\big)$
by $\phi_{\sigma,\tau}$, where $\sigma$ is a~$3$-cycle and
$\tau$ is a transposition such that $\sigma \tau = g$. Then $\phi_{\sigma,\tau} \equiv 0$.
\end{Lemma}

\begin{proof}
The proof proceeds by considering the same exhaustive cases as in Lemma~\ref{lemma:simplification*2L2},
but using the definition of $\kappa^C_{\tri}$ to show in fact $\phi_{\sigma,\tau}\equiv 0$ in cases~(2) and~(3).
Showing that $\phi^{*}_{\sigma,\tau}(u,v,w)=0$ when
$u,v,w \in W$, $u,v,w \in W^*$, $u,v\in V^{\tau}$, or $u \in V^{\tau}\cap V^{\sigma}$
proceeds exactly as~in~the proof of Lemma~\ref{lemma:simplification*2L2}
since the methods did not depend on anything about $\kappa^*_\sigma$ other than
$V^{\sigma} \subseteq \ker \kappa^*_{\sigma}$.

As in the proof of case~(2) in Lemma~\ref{lemma:simplification*2L2}, assume
$u \in V^{\tau} \setminus V^{\sigma}$ and $v \notin V^\tau$ and note
	\begin{gather*}
	\phi_{\sigma,\tau}(v,\actbylessspace{\tau}v,u^*)=0,
	\end{gather*}
and
	\begin{gather*}
	\phi_{\sigma,\tau}(u,v,v^*)=2\kappa^C_{\sigma}\big(u,\kappa^L_{\tau}(v,v^*)\big).
	\end{gather*}
Using bilinearity and $V^{\sigma} \subseteq \ker \kappa^*_{\sigma}$, the right hand
side is a linear combination of expressions $\kappa^*_{\sigma}\big(u,\actby{h}u\big)$ for
$h \in \langle \sigma \rangle$. The appropriate coefficients in terms of $a$, $a^\perp$, $b$, $b^\perp$
can be described in~terms of the indicator function for the fixed space of $\tau$
and depend on whether $u \in W^*$ or $u \in W$.
Also, $\sum_{h \in \langle \sigma \rangle} \actby{h}u \in V^{\sigma} \subseteq \ker \kappa^*_{\sigma}$.
Thus for $u \in W^*$,
\begin{gather*}
\kappa^C_{\sigma}\big(u,\kappa^L_{\tau}(v,v^*)\big)
=\big(a^\perp-a\big)\sum\limits_{h \in \langle \sigma \rangle}
\delta_{\tau}\big(\actby{h}u\big)\kappa^C_{\sigma}\big(u,\actby{h}u\big)
+\big(b^\perp-b\big)\sum\limits_{h \in \langle \sigma \rangle}
\delta_{\tau}\big(\actby{h}u^*\big)\kappa^C_{\sigma}\big(u,\actby{h}u^*\big)
\\ \hphantom{\kappa^C_{\sigma}\big(u,\kappa^L_{\tau}(v,v^*)\big)}
=\big(a^\perp-a\big)\big[\delta_{\tau}(\actbylessspace{\sigma}u)\kappa^C_{\sigma}(u,\actbylessspace{\sigma}u)
+\delta_{\tau}\big(\actbylessspace{\sigma^2}u\big)\kappa^C_{\sigma}\big(u,\actbylessspace{\sigma^2}u\big)\big]
\\ \hphantom{\kappa^C_{\sigma}\big(u,\kappa^L_{\tau}(v,v^*)\big)=}
\phantom{,}+\big(b^\perp-b\big)\big[\delta_{\tau}(\actbylessspace{\sigma}u)\kappa^C_{\sigma}(u,\actbylessspace{\sigma}u^*)
+\delta_{\tau}\big(\actbylessspace{\sigma^2}u\big)\kappa^C_{\sigma}\big(u,\actbylessspace{\sigma^2}u^*\big)\big]
\\ \hphantom{\kappa^C_{\sigma}\big(u,\kappa^L_{\tau}(v,v^*)\big)}
=\big[\delta_{\tau}(\actbylessspace{\sigma}u)-\delta_{\tau}\big(\actbylessspace{\sigma^2}u\big)\big]
\big[\big(a^\perp\!-a\big)\big(b^\perp\!-b\big)^2-\big(b^\perp\!-b\big)\big(a^\perp\!-a\big)\big(b^\perp\!-b\big) \big] =0,
\end{gather*}
and hence $\phi_{\sigma,\tau}(u,v,v^*)=0$.
A similar calculation shows that
	\begin{gather*}
\phi_{\sigma,\tau}(u,v,\actbylessspace{\tau}v^*)=-\phi_{\sigma,\tau}(u,v,v^*)=0
\end{gather*}
and applying $G$-invariance
leads to the same conclusions when $u \in W$.
Then~\eqref{tauinv} also implies
\begin{gather*}
\phi_{\sigma,\tau}(u,\actby{\tau}v,\actby{\tau}v^*)=-\phi_{\sigma,\tau}(u,\actby{\tau}v,v^*) =\phi_{\sigma,\tau}(u,v,v^*)=0.
\end{gather*}
	
For case~(3) assume $v \notin V^{\tau}$ and note
that as in the proof of Lemma~\ref{lemma:simplification*2L2},
\begin{gather*}
\phi_{\sigma,\tau}(v,\actbylessspace{\tau}v,v^*)
=-2\kappa^C_{\sigma}\big(v+\actbylessspace{\tau}v,\kappa^L_{\tau}(v,v^*)\big).
\end{gather*}
Since $\kappa^C_{\sigma}\big(\actbylessspace{\tau}v,\kappa^L_{\tau}(v,v^*)\big)
=\kappa^C_{\sigma}\big(\actbylessspace{\tau}v,\kappa^L_{\tau}(\actbylessspace{\tau}v, \actbylessspace{\tau}v^*)\big)
=\kappa^C_{\sigma}\big(v,\kappa^L_{\tau}(v,v^*)\big)
=0$ by the same calculation as in
case~(2) except with $\actbylessspace{\tau}v$ and $v$ in place of $u$, it follows that
\begin{gather*}
\phi_{\sigma,\tau}(v,\actbylessspace{\tau}v,v^*)=
-2\big[\kappa^C_{\sigma}\big(v,\kappa^L_{\tau}(v,v^*)\big)+
\kappa^C_{\sigma}\big(\actbylessspace{\tau}v,\kappa^L_{\tau}(v,v^*)\big)\big]=0.
\end{gather*}
Lastly,~\eqref{tauinv} yields
\begin{gather*}
\phi_{\sigma,\tau}(v,\actby{\tau}v,\actby{\tau}v^*)=-\phi_{\sigma,\tau}(v,\actby{\tau}v,v^*)=0.\tag*{\qed}
\end{gather*}
\renewcommand{\qed}{}
\end{proof}

The case where $*=L$ is included in the preliminary calculations of the following
lemma because it will be useful as a starting point in the proof of
Theorem~\ref{thm:DOAMapsCombined}. We~note that when $n=2$ conditions~\eqref{Obstr2PhiC1C2C3L2=0a} and~\eqref{Obstr2PhiC1C2C3L2=0b} need to be modified to $a(\alpha+\beta)=b(\alpha+\beta)=0$.

\begin{Lemma}\label{lemma:simplification*1L2}
Let $\kappa^C_1$ and $\kappa^L_{\refl}$ be as in
Definitions~$\ref{def:kappa1}$ and~$\ref{def:kapparefl}$,
with parameters $\alpha, \beta \in \mathbb{C}$ and $a, a^\perp, b, b^\perp \in \mathbb{C}$
respectively.
Denote a term of the component $\phi_g$ of $\phi\big(\kappa^*_1,\kappa^L_{\refl}\big)$
by $\phi^*_{1,g}$, where $*=C$ or $*=L$ and $g$ is a transposition.
Then $\phi^C_{1,g}=0$ if and only if the following conditions hold
\begin{alignat}{3}
\label{Obstr2PhiC1C2C3L2=0a}
&\alpha a + \beta\big(a+(n-2)a^\perp\big)=0, \qquad &&\alpha a^\perp + \beta\big(2a+(n-3)a^\perp\big)=0,&
\\
\label{Obstr2PhiC1C2C3L2=0b}
&\alpha b + \beta\big(b+(n-2)b^\perp\big)=0, \qquad &&\alpha b^\perp + \beta\big(2b+(n-3)b^\perp\big)=0.&
\end{alignat}
\end{Lemma}

\begin{proof}
As in~Section~\ref{InvLieBrackets}, it suffices to compute $\phi^*_{1,g}$ on basis triples
of the following forms for $1 \leq i, j, k \leq n$.
\begin{alignat*}{3}
&\text{1.\ All basis vectors in $W$ or in $W^*$ and $i$, $j$, $k$ distinct:} &\hspace{10pt}
& (x_i,x_j,x_k), & &\quad (y_i,y_j,y_k). \\
&\text{2.\ Two basis vectors in $W$ or in $W^*$ and $i$, $j$, $k$ distinct:} &\hspace{10pt}
& (x_i,x_j,y_k), & &\quad (y_i,y_j,x_k). \\
&\text{3.\ Two basis vectors in $W$ or $W^*$ and $i$, $j$ distinct:} &\hspace{10pt}
& (x_i,x_j,y_j), & &\quad (y_i,y_j,x_j).
\end{alignat*}

\textit{Case}~1.
For any distinct indices $i$, $j$, $k$ with $1 \leq i, j, k \leq n$, using~\eqref{phi} and
Definition~\ref{def:kapparefl} of~$\kappa^L_{\refl}$ it is immediate
that $\phi^*_{1,g}(x_i,x_j,x_k)=0$, and in similar fashion $\phi^*_{1,g}(y_i,y_j,y_k)=0$,
for any (distinct $i$, $j$, $k$ with) $1 \leq i, j, k \leq n$.
Thus this case imposes no conditions on any parameters.

\textit{Case}~2.
For any distinct indices $i$, $j$, $k$ with $1 \leq i, j, k \leq n$,
using the definitions of~$\kappa^L_{\refl}$ and~$\kappa^C_1$, bilinearity,
and skew-symmetry yields
\begin{align*}
\phi^*_{1,g}(x_i,x_j,y_k)
&=\begin{cases}
2\kappa^*_1\big(x_j,ax_{ik}+a^\perp x_{ik}^\perp + b y_{ik} + b^\perp y_{ik}^\perp\big) &\text{if}\quad g=(ik),
\\
-2\kappa^*_1\big(x_i,ax_{jk}+a^\perp x_{jk}^\perp + b y_{jk} + b^\perp y_{jk}^\perp\big) &\text{if}\quad g=(jk), \\
0 & \text{otherwise}.
\end{cases}
\\
&=\begin{cases}
2\big[\alpha b^\perp + \beta\big(2b+(n-3)b^\perp\big)\big] & \text{if}\quad g=(ik)\quad \text{and}\quad *=C, \\
-2\big[\alpha b^\perp + \beta\big(2b+(n-3)b^\perp\big)\big] &\text{if}\quad g=(jk) \quad \text{and}\quad *=C, \\
 0 & \text{otherwise}.
\end{cases}
\end{align*}

Interchanging the roles of $x$ and $y$ and recomputing yields that for any
distinct indices $i$, $j$, $k$ with $1 \leq i, j, k \leq n$,
\begin{gather*}
\phi^*_{1,g}(y_i,y_j,x_k)=
\begin{cases}
2\big[\alpha a^\perp + \beta\big(2a+(n-3)a^\perp\big)\big] & \text{if}\quad g=(ik) \quad\text{and}\quad *=C, \\
	 -2\big[\alpha a^\perp + \beta\big(2a+(n-3)a^\perp\big)\big] & \text{if}\quad g=(jk) \quad\text{and}\quad *=C, \\
 0 & \text{otherwise}.
 \end{cases}
\end{gather*}

\textit{Case}~3.
For any distinct indices $i$, $j$ with $1 \leq i, j \leq n$,
using the definitions of~$\kappa^L_{\refl}$ and~$\kappa^C_1$, bilinearity,
and skew-symmetry yields
\begin{align*}
\phi^*_{1,g}(x_i,x_j,y_j)
&=\begin{cases}
2\kappa^*_1\big(x_i+x_j,ax_{ij}+a^\perp x_{ij}^\perp + by_{ij} + b^\perp y_{ij}^\perp\big) &\text{if}\quad g=(ij), \\
2\kappa^*_1\big(x_i,ax_{jk}+a^\perp x_{jk}^\perp + b y_{jk} + b^\perp y_{jk}^\perp\big) &\text{if}\quad g=(jk), \\
 0 & \text{otherwise}.
\end{cases}
\\
&=\begin{cases}
4\big[\alpha b + \beta\big(b+(n-2)b^\perp\big)\big] & \text{if}\quad g=(ij) \quad\text{and}\quad *=C, \\
2\big[\alpha b^\perp + \beta\big(2b+(n-3)b^\perp\big)\big] & \text{if}\quad g=(jk) \quad\text{and}\quad *=C,
\\
0 & \text{otherwise}.
\end{cases}
\end{align*}

Interchanging the roles of $x$ and $y$ and recomputing yields that for any
distinct indices $i$, $j$ with $1 \leq i, j \leq n$,
\begin{gather*}
\phi^*_{1,g}(y_i,y_j,x_j)=
\begin{cases}
4\big[\alpha a + \beta\big(a+(n-2)a^\perp\big)\big] & \text{if}\quad g=(ij) \quad\text{and}\quad *=C, \\
2\big[\alpha a^\perp + \beta\big(2a+(n-3)a^\perp\big)\big] & \text{if}\quad g=(jk) \quad\text{and}\quad *=C, \\
 0 & \text{otherwise}.
 \end{cases}
\end{gather*}

Setting the results in cases~2 and~3 equal to zero yields
conditions~\eqref{Obstr2PhiC1C2C3L2=0a} and~\eqref{Obstr2PhiC1C2C3L2=0b}.
\end{proof}

\begin{Corollary}\label{cor:phiC3L2=0}
 Let $\kappa^L_{\refl}$ and $\kappa^C_{\tri}$ be as in Definitions~$\ref{def:kapparefl}$
 and~$\ref{def:kappatri}$, with common parameters $a,a^\perp,b,b^\perp \in\mathbb{C}$. For~every $g\in S_n$, the component $\phi_g$ of $\phi\big(\kappa^C_{\refl}+\kappa^C_{\tri},\kappa^L_{\refl}\big)$
 is identically zero. The component $\phi_g$ of $\phi\big(\kappa^C_1,\kappa^L_{\refl}\big)$
 is zero for all $g \in S_n$ if and only if conditions~\eqref{Obstr2PhiC1C2C3L2=0a}
 and~\eqref{Obstr2PhiC1C2C3L2=0b} given in Lemma~$\ref{lemma:simplification*1L2}$ hold.
\end{Corollary}

\begin{proof}
In fact, each component of $\phi\big(\kappa^C_{\refl},\kappa^L_{\refl}\big)$
is identically zero by Propositions~\ref{prop:zerocases} and~\ref{prop:3-cycle},
and each component of $\phi\big(\kappa^C_{\tri},\kappa^L_{\refl}\big)$ is identically zero
by Lemma~\ref{lemma:simplificationC3L2}. By~Lemma~\ref{lemma:simplification*1L2}
we have $\phi\big(\kappa^C_1,\kappa^L_{\refl}\big)=0$ if and only if
conditions~\eqref{Obstr2PhiC1C2C3L2=0a} and~\eqref{Obstr2PhiC1C2C3L2=0b}
are satisfied because the component $\phi_g$ is identically zero for $g$ not a
transposition by the definitions of $\kappa^C_1$ and $\kappa^L_{\refl}$.
\end{proof}

Conditions~\eqref{Obstr2PhiC1C2C3L2=0a} and~\eqref{Obstr2PhiC1C2C3L2=0b}
given in Lemma~\ref{lemma:simplification*1L2}
give rise to a variety that controls the para\-me\-ter space for the
family of maps in Theorem~\ref{thm:DOAMapsDoubledPerm}
and Drinfeld orbifold algebras in Theorem~\ref{permexamples}. We~conjecture that this projective variety has dimension four based on computations done
for a few specific values of $n$ in Macaulay2~\cite{M2} with the graded reverse
lexicographic monomial ordering and the
parameter order $a$, $a^\perp$, $b$, $b^\perp$, $\alpha$, $\beta$, $c$.

\subsection{Drinfeld orbifold algebra maps}
Now we use the details of clearing the obstructions from
Section~\ref{ObstructionDetails} to describe all Drinfeld orbi\-fold
algebra maps with linear part supported only off the identity.
The corresponding Drinfeld orbifold algebras are given in Theorem~\ref{permexamples}.

\begin{Theorem}\label{thm:DOAMapsDoubledPerm}
For $S_n$ $(n \geq 3)$ acting on $V=W^* \oplus W \cong \mathbb{C}^{2n}$ by the
doubled permutation representation, the Drinfeld orbifold algebra maps
supported only off the identity are precisely the maps of the form
$\kappa=\kappa^L_{\refl}+\kappa^C_{\tri}+\kappa^C_{\refl}+\kappa^C_1$, with
$\kappa^L_{\refl}$ and $\kappa^C_{\refl}$ as in Definition~$\ref{def:kapparefl}$,
$\kappa^C_{\tri}$ as in Definition~$\ref{def:kappatri}$, $\kappa^C_1$
as in Definition~$\ref{def:kappa1}$, and with the parameters
$a$, $a^\perp$, $b$, $b^\perp$, $c$, $\alpha$, and $\beta$ satisfying
these conditions derived in Lemma~$\ref{lemma:simplification*1L2}$:
\begin{enumerate}	\itemsep=0pt
	\item[\eqref{Obstr2PhiC1C2C3L2=0a}] \ $\alpha a + \beta\big(a+(n-2)a^\perp\big)=0, \qquad \alpha a^\perp + \beta\big(2a+(n-3)a^\perp\big)=0$,
	\item[\eqref{Obstr2PhiC1C2C3L2=0b}] \ $\alpha b + \beta\big(b+(n-2)b^\perp\big)=0, \qquad \alpha b^\perp + \beta\big(2b+(n-3)b^\perp\big)=0$.
\end{enumerate}
In particular, $\kappa=\kappa^L_{\refl}+\kappa^C_{\tri}+\kappa^C_{\refl}$ is
always a Drinfeld orbifold algebra map.
\end{Theorem}

\begin{proof}
Suppose $\kappa^L$ is a pre-Drinfeld orbifold algebra map supported only off
the identity. By~Corollary~\ref{preDOAmaps} we must have $\kappa^L=\kappa^L_{\refl}$
for some parameters $a,a^\perp,b,b^\perp \in \mathbb{C}$ as in Definition~\ref{def:kapparefl}.
It remains to find all $G$-invariant maps $\kappa^C$ such that
properties~\eqref{PBWconditions-iii} and~\eqref{PBWconditions-iv}
of a Drinfeld orbifold algebra map also hold.

First we find a particular lift.
\begin{itemize}\itemsep=0pt
	\item {\it First obstruction.} Propositions~\ref{prop:zerocases} and~\ref{prop:3-cycle} give the value
	of $\phi\big(\kappa^L_{\refl},\kappa^L_{\refl}\big)$. These values sug\-gest how to construct
	the $S_n$-invariant map $\kappa^C_{\tri}$ such that property~\eqref{PBWconditions-iii}
	holds, as~given in Definition~\ref{def:kappatri}. Proposition~\ref{prop:phiL2L2=2psiC3}
	then verifies that $\phi\big(\kappa^L_{\refl},\kappa^L_{\refl}\big)$ and
	$2\psi\big(\kappa^C_{\tri}\big)$ are indeed equal.
	\item {\it Second obstruction.} By Lemma~\ref{lemma:simplificationC3L2}, we have
	that $\phi\big(\kappa^C_{\tri},\kappa^L_{\refl}\big)=0$.
\end{itemize}
	Thus $\kappa=\kappa^L_{\refl}+\kappa^C_{\tri}$ is 	a Drinfeld orbifold algebra map.

Next, we modify this particular lift to obtain all possible lifts.
Let $\kappa^C$ be any $G$-invariant constant $2$-cochain.
\begin{itemize}\itemsep=0pt
 \item {\it First obstruction. } Since $\phi\big(\kappa^L_{\refl},\kappa^L_{\refl})=2\psi(\kappa^C_{\tri}\big)$,
 it follows that $\phi\big(\kappa^L_{\refl},\kappa^L_{\refl}\big)=2\psi\big(\kappa^C\big)$
 if and only if $\psi\big(\kappa^C-\kappa^C_{\tri}\big)=0$. By~Corollary~\ref{mixedJacobi} this occurs
 if and only if $\kappa^C-\kappa^C_{\tri}=\kappa^C_{\refl}+\kappa^C_1$, with $\kappa^C_{\refl}$
 as in Definition~\ref{def:kapparefl} for some parameter $c\in\mathbb{C}$ and $\kappa^C_1$
 as in Definition~\ref{def:kappa1} for some parameters $\alpha,\beta\in\mathbb{C}$.
 \item {\it Second obstruction. }
 By Corollary~\ref{cor:phiC3L2=0},
 $\phi\big(\kappa^C_1+\kappa^C_{\refl}+\kappa^C_{\tri},\kappa^L_{\refl}\big)=0$
 if and only if~\eqref{Obstr2PhiC1C2C3L2=0a} and~\eqref{Obstr2PhiC1C2C3L2=0b} are satisfied.
 This occurs in particular when $\kappa^C_1 \equiv 0$ in which case
 $\kappa^C_{\refl}+\kappa^C_{\tri}$ clears the second obstruction and
 $\kappa^L_{\refl}$ lifts with no restrictions on the parameters $a$,
 $a^\perp$, $b$, $b^\perp$, or $c$.
 \end{itemize}
 Thus the lifts of $\kappa^L_{\refl}$ to a Drinfeld orbifold algebra map are
 precisely the maps of the form
 $\kappa=\kappa^L_{\refl}+\kappa^C_1+\kappa^C_{\refl}+\kappa^C_{\tri}$ satisfying
 conditions~\eqref{Obstr2PhiC1C2C3L2=0a} and~\eqref{Obstr2PhiC1C2C3L2=0b}.
\end{proof}

Lastly, by specifying parameters we obtain some Drinfeld orbifold algebra
maps that are supported both on and off the identity. The proof uses results
related to clearing obstructions that appeared in~Section~\ref{InvLieBrackets}
and Section~\ref{ObstructionDetails}.

\begin{Theorem}
\label{thm:DOAMapsCombined}
For $S_n$ $(n \geq 3)$ acting on $V=W^* \oplus W \cong \mathbb{C}^{2n}$ by the
doubled permutation representation, there are Drinfeld orbifold algebra maps
of the form $\kappa=\kappa^L+\kappa^C$
with linear part $\kappa^L=\kappa^L_1+\kappa^L_{\refl}$ and
constant part $\kappa^C=\kappa^C_1+\kappa^C_{\refl}+\kappa^C_{\tri}$, where
\begin{enumerate}\itemsep=0pt
\item[$(1)$] $\kappa^L_1$ is as described in Definition~$\ref{def:kappa1}$ with $a_i=b_i=0$
for $i=1, 2, 3, 5, 6$ and $a_4=a_7$ and~$b_4=b_7$ are not both zero,
\item[$(2)$] $\kappa^C_1$ is as described in Definition~$\ref{def:kappa1}$ with $\alpha=\beta$,
\item[$(3)$] $\kappa^L_{\refl}$ and $\kappa^C_{\refl}$ are as in Definition~$\ref{def:kapparefl}$ and
$\kappa^C_{\tri}$ is as in Definition~$\ref{def:kappatri}$ with
$2a+(n-2)a^\perp=2b+(n-2)b^\perp=0$, but $a$, $a^\perp$, $b$, and $b^\perp$ not all zero.
\end{enumerate}
\end{Theorem}

\begin{proof}
As in the proof of Theorem~\ref{thm:DOAMapsDoubledPerm},
even without the given parameter choices we have
$\phi\big(\kappa^L_{\refl},\kappa^L_{\refl}\big)=2\psi\big(\kappa^C_{\tri}\big)=
2\psi\big(\kappa^C_1+\kappa^C_{\refl}+\kappa^C_{\tri}\big)$.
Setting $a_i=b_i=0$ for $i=1, 2, 3, 5, 6$, $a_4=a_7$, and $b_4=b_7$
in~\eqref{ClearFirstObstr} and in the values of $\phi^L_{g,1}(u,v,w)$
in~Sections~\ref{AllVecInW}--\ref{TwoVecInWTwoIndices} shows
$\phi\big(\kappa^L_1,\kappa^L_1\big)=0$ and $\phi\big(\kappa^L_{\refl},\kappa^L_1\big)=0$ respectively.
Using the forms of $\phi^L_{1,g}$ given in Lemma~\ref{lemma:simplification*1L2} and the assumptions that $a_4=a_7$ and $b_4=b_7$ yield that
\begin{gather*}
\phi^L_{1,g}(x_i,x_j,y_k)=
\begin{cases}
2\big(2b + (n-2)b^\perp\big)\big(a_4x_{[n]}+b_4y_{[n]}\big) & \text{if}\quad g=(ik), \\
-2\big(2b + (n-2)b^\perp\big)\big(a_4x_{[n]}+b_4y_{[n]}\big) & \text{if}\quad g=(jk), \\
0 & \text{otherwise},
\end{cases}
\\
\phi^L_{1,g}(x_i,x_j,y_j)=
\begin{cases}
4\big(2b + (n-2)b^\perp\big)\big(a_4x_{[n]}+b_4y_{[n]}\big) & \text{if}\quad g=(ij), \\
2\big(2b + (n-2)b^\perp\big)\big(a_4x_{[n]}+b_4y_{[n]}\big) & \text{if}\quad g=(jk), \\
 0 & \text{otherwise},
\end{cases}
\end{gather*}
and
similarly for $\phi^L_{1,g}(y_i,y_j,x_k)$ and $\phi^L_{1,g}(y_i,y_j,x_j)$,
but replacing $b$ with $a$ and $b^\perp$ with $a^\perp$. By~the hypothesis on $a$, $a^\perp$, $b$, and $b^\perp$, all of these are zero,
so $\phi\big(\kappa^L_1,\kappa^L_{\refl}\big)=0$
and hence
\begin{gather*}
\phi\big(\kappa^L_1+\kappa^L_{\refl},\kappa^L_1+\kappa^L_{\refl}\big)=
\phi\big(\kappa^L_{\refl},\kappa^L_{\refl}\big)=
2\psi\big(\kappa^C_{\tri}\big)=
2\psi\big(\kappa^C_1+\kappa^C_{\refl}+\kappa^C_{\tri}\big).
\end{gather*}
Thus with the given parameter choices
$\kappa^C_1+\kappa^C_{\refl}+\kappa^C_{\tri}$
clears the first obstruction for $\kappa^L_1+\kappa^L_{\refl}$.

We claim $\kappa^C_1+\kappa^C_{\refl}+\kappa^C_{\tri}$ also clears the second obstruction
for $\kappa^L_1+\kappa^L_{\refl}$ because
\begin{gather*}
\phi\big(\kappa^C_1+\kappa^C_{\refl}+\kappa^C_{\tri},\kappa^L_1+\kappa^L_{\refl}\big)=0.
\end{gather*}
By the assumptions on $a_i$ and $b_i$ for $1\leq i \leq 7$
conditions~\eqref{ClearSecondObstr1C} and~\eqref{ClearSecondObstr1Ref} are satisfied and thus
$\phi\big(\kappa^C_1+\kappa^C_{\refl},\kappa^L_1\big)=0$.
Also $\phi\big(\kappa^C_{\tri},\kappa^L_1\big)=0$. This is because
for any $1\leq i,j,k \leq n$ and any $3$-cycle $g$,
by $a_1=a_2=b_1=b_2=0$ we have
\begin{gather*}
\phi^C_{g,1}(x_i,x_j,x_k)=\phi^C_{g,1}(y_i,y_j,y_k)=0,
\end{gather*}
and by $a_4=a_7$ and $b_4=b_7$ we have
\begin{align*}
\phi^C_{g,1}(x_i,x_j,y_k)
=\kappa^C_g\big(x_i,\kappa^L_1(x_j,y_k)\big)+\kappa^C_g\big(x_j,\kappa^L_1(y_k,x_i)\big)
=\kappa^C_g(x_i-x_j,a_4 x_{[n]} + b_4 y_{[n]}) =0
\end{align*}
because $\kappa^C_g\big(x_i,\actby{g}x_i+\actbylessspace{g^{-1}}x_i\big)=0$
regardless of whether $i \in V^g$ by Definition~\ref{def:kappatri},
and similarly $\phi^C_{g,1}(y_i,y_j,x_k)=0$. By~Corollary~\ref{cor:phiC3L2=0} we know that
$\phi\big(\kappa^C_{\refl}+\kappa^C_{\tri},\kappa^L_{\refl}\big)=0$ in general,
and that $\phi\big(\kappa^C_1,\kappa^L_{\refl}\big)=0$ because
conditions~\eqref{Obstr2PhiC1C2C3L2=0a} and~\eqref{Obstr2PhiC1C2C3L2=0b}
are satisfied when
$\alpha=\beta$ and $2a+(n-2)a^\perp=2b+(n-2)b^\perp=0$.

Thus with the given choices of parameters,
$\kappa^C_1+\kappa^C_{\refl}+\kappa^C_{\tri}$ clears the second obstruction
as well and lifts $\kappa^L_1+\kappa^L_{\refl}$ to a Drinfeld orbifold algebra map.
\end{proof}

\section[Drinfeld orbifold algebra maps that deform S(h* + h)S_n]
{Drinfeld orbifold algebra maps that deform $\boldsymbol{S(\mathfrak{h}^*\oplus\mathfrak{h})\# S_n}$}\label{DoubledStandardDefs}

We now use the results in Sections~\ref{InvLieBrackets} and~\ref{Lifting}
on Lie and Drinfeld orbifold algebra
maps that produce deformations of the skew group algebra $S(W^* \oplus W)\# S_n$
in order to understand which maps produce deformations of $S(\mathfrak{h}^*\oplus\mathfrak{h})\# S_n$.

In contrast to the complicated families of Lie orbifold algebras and
maps in Theorems~\ref{thm:LOAMapsDoubledPerm} and~\ref{Lieorbifold},
when $S_n$ instead acts on its doubled standard subrepresentation
$\mathfrak{h}^*\oplus\mathfrak{h}$ there are no Lie orbifold algebra maps with nonzero linear
part (Theorem~\ref{thm:LOAMapsDoubledStd}).
However, Theorem~\ref{thm:DOAMapsDoubledStd} describes a~three-parameter
family of~Drinfeld orbifold algebra maps that do provide polynomial degree
one deformations generalizing the $\mathfrak{sl}_n$-type rational Cherednik
algebras $H_{0,c}$ (see also Theorem~\ref{refexamples}). We~begin with a result that applies to any finite group and
provides conditions under which we can combine Drinfeld orbifold algebra maps for subrepresentations into a map for their direct~sum.

\begin{Proposition}	\label{prop:directsum}
Let $G$ be a finite group acting linearly on finite-dimensional vector spaces $U_1, U_2,\ldots,U_r$.
Given Drinfeld orbifold algebra maps $\kappa|_{U_i}$ for $G$ acting on $U_i$ $(i=1,\ldots,r)$, define $\kappa$ on $\bigwedge^2 \big(\bigoplus_{i=1}^r U_i\big)$ so that $\kappa$ agrees with $\kappa|_{U_i}$ for pairs of vectors from the same $U_i$ and is zero on mixed pairs, i.e., $\kappa(U_i,U_j)=0$ for $i\neq j$. Then $\kappa$ is a Drinfeld orbifold algebra map for $G$ acting on $\bigoplus_{i=1}^r U_i$ if and only if whenever $i\neq j$ all group elements in the support of $\kappa|_{U_i}$ act trivially on $U_j$.
\end{Proposition}

\begin{proof}
 For $i=1,\ldots,r$, suppose $\kappa|_{U_i}$ is a Drinfeld orbifold algebra map for $G$ acting on $U_i$ and define $\kappa$ on $\bigwedge^2\big(\bigoplus_{i=1}^r U_i\big)$ as above.
 Conditions~\eqref{PBWconditions-0} and~\eqref{PBWconditions-i} of the definition of a Drinfeld orbifold algebra map are straightforward to verify. For~conditions~\eqref{PBWconditions-ii}--\eqref{PBWconditions-iv}, consider a triple of vectors from $\bigoplus_{i=1}^r U_i$. If~all three vectors are from the same $U_i$, then equations~\eqref{PBWconditions-ii}--\eqref{PBWconditions-iv} hold by virtue of $\kappa|_{U_i}$ being a Drinfeld orbifold algebra map. If~the three vectors are from~$U_i$,~$U_j$, and~$U_k$ with $i$, $j$, $k$ distinct, then conditions~\eqref{PBWconditions-ii}--\eqref{PBWconditions-iv} are easily seen to hold because~$\kappa$ is defined to be zero on pairs of vectors from different summands.

 By multilinearity and skew-symmetry, all that remains is to examine the case where two vectors, say $u$ and $v$, are from the same $U_i$ and the third vector, say $w$, is from some $U_j$ with $j\neq i$. Recall that
 \begin{gather*}
 \phi_{x,y}^*(u,v,w)=\kappa_x^*\big(u+\actby{y}u,\kappa_y^L(v,w)\big) +\kappa_x^*\big(v+\actby{y}v,\kappa_y^L(w,u)\big)+\kappa_x^*\big(w+\actby{y}w,\kappa_y^L(u,v)\big).
 \end{gather*}
 In the present case, the first two terms are zero because $\kappa^L_y(U_i,U_j)=0$, and the last term is zero because $\kappa_y^L(U_i,U_i)\subseteq U_i$ and $\kappa^*_x(U_j,U_i)=0$. Thus $\phi_{x,y}^*(u,v,w)=0$, which implies equation~\eqref{PBWconditions-iv} is satisfied for all $g$ in~$G$. We~also see~\eqref{PBWconditions-ii} and~\eqref{PBWconditions-iii} will be satisfied if and only if for all $g\in G$, we have $\psi_g^*(u,v,w)=0$ for $*=L$ and $*=C$, respectively.

 To this end, recall that
 \begin{gather*}
 \psi_g^*(u,v,w)=\kappa_g^*(u,v)\big(\actby{g}w-w\big)+\kappa_g^*(v,w)\big(\actby{g}u-u\big)+ \kappa_g^*(w,u)\big(\actby{g}v-v\big).
 \end{gather*}
 Continuing with $u,v\in U_i$ and $w\in U_j$, we see that $\psi_g^*(u,v,w)=\kappa_g^*(u,v)(\actby{g}w-w)$ because $\kappa_g^*(U_i,U_j)=0$.
 Thus for conditions~\eqref{PBWconditions-ii} and~\eqref{PBWconditions-iii} to hold for all $g$ in~$G$ and all triples of this type, it is both necessary and sufficient for group elements in the support of $\kappa|_{U_i}$ to act trivially on $U_j$ whenever $i\neq j$.
\end{proof}

 As a corollary, in the case of two summands $U_1$ and $U_2$ with $G$ acting trivially on $U_1$ and $\kappa|_{U_1}\equiv 0$ (so that the support of $\kappa|_{U_1}$ is empty), we have:

\begin{Corollary}\label{prop:Extension}
	Let $G$ be a finite group acting linearly on a vector space $V=U_1\oplus U_2$, where each $U_i$ is a subrepresentation. If~$G$ acts trivially on $U_1$, then every Drinfeld orbifold algebra map $\kappa|_{U_2}$ on $\bigwedge^2 U_2$ extends to a Drinfeld orbifold algebra map $\kappa$ on $\bigwedge^2 V$ with $\im\kappa^L\subseteq U_2$ and such that $U_1\subseteq \ker\kappa$.
\end{Corollary}

We now use Corollary~\ref{prop:Extension} to show there are no Drinfeld orbifold algebra maps with linear part
supported only on the identity for $S_n$ acting on $\mathfrak{h}^* \oplus \mathfrak{h}$.

\begin{Theorem}\label{thm:LOAMapsDoubledStd}
For $S_n$ $(n \geq 3)$ acting on $\mathfrak{h}^*\oplus\mathfrak{h} \cong\mathbb{C}^{2n-2}$ by the doubled
standard representation there are no degree-one Lie orbifold algebra maps.
\end{Theorem}

\begin{proof}
If there were a Lie orbifold algebra map $\kappa^L+\kappa^C$ for the
doubled standard representation with $\kappa^C\colon \bigwedge^2 (\mathfrak{h}^*\oplus\mathfrak{h}) \to \mathbb{C} S_n$
and nonzero $\kappa^L\colon \bigwedge^2 (\mathfrak{h}^*\oplus\mathfrak{h}) \to \mathfrak{h}^*\oplus\mathfrak{h}$, then
it could be extended as in Corollary~\ref{prop:Extension} to yield a
Lie orbifold algebra map for $S_n$ acting on $V=W^* \oplus W \cong \mathbb{C}^{2n}$ via
the doubled permutation representation. The possible forms of Lie orbifold algebra maps
$\kappa$ for the doubled permutation representation are controlled by Theorem~\ref{thm:LOAMapsDoubledPerm}, which includes the PBW conditions~$\gamma_1=\gamma_2=\gamma_4=\gamma_5=0$ in~\eqref{ClearFirstObstr}. We~will use these to show that in fact $\kappa^L\equiv 0$
by first imposing the image constraint $\im \kappa^L \subseteq \mathfrak{h}^*\oplus\mathfrak{h}$ and the kernel constraint $\iota^*\oplus\iota\subseteq\ker\kappa$
from Corollary~\ref{prop:Extension}.

First, use
$x_i=\bar{x}_i+\tfrac{1}{n}x_{[n]}$ and $y_i=\bar{y}_i+\tfrac{1}{n}y_{[n]}$
in~\eqref{kappaLxx}--\eqref{kappaLxiyj} to write the values of $\kappa^L_1$
according to the decomposition $V\cong\mathfrak{h}^* \oplus \iota^* \oplus \mathfrak{h} \oplus \iota$:
\begin{gather*}
\kappa^L_1(x_i,x_j)=a_1(\bar{x}_i-\bar{x}_j)+b_1(\bar{y}_i-\bar{y}_j), \\
\kappa^L_1(y_i,y_j)=a_2(\bar{x}_i-\bar{x}_j)+b_2(\bar{y}_i-\bar{y}_j), \\
\kappa^L_1(x_i,y_i)=a_3\bar{x}_i+\frac{1}{n}(a_3+na_4)\n{x}+b_3 \bar{y}_i+\frac{1}{n}(b_3+nb_4)\n{y}, \\
\kappa^L_1(x_i,y_j)=a_5 \bar{x}_i+a_6 \bar{x}_j+\frac{1}{n}(a_5+a_6+na_7)\n{x}+b_5 \bar{y}_i+b_6 \bar{y}_j+\frac{1}{n}(b_5+b_6+nb_7)\n{y}.
\end{gather*}
Thus $\im \kappa^L \subseteq \mathfrak{h}^*\oplus\mathfrak{h}$ implies
\begin{alignat*}{3}
& a_3+na_4=0, \qquad && b_3+nb_4=0,& \\
&a_5+a_6+na_7=0, \qquad && b_5+b_6+nb_7=0.&
\end{alignat*}

Second, impose the extension conditions $\kappa^L(w_i,v_j)=\kappa^L(v_i,v_j)=0$ for
$w_i$ in the basis $\{x_{i+1}-x_i, y_{i+1}-y_i \mid 1 \leq i \leq n-1 \}$
of $\mathfrak{h}^*\oplus\mathfrak{h}$ and $v_i, v_j$ in the basis $\{x_{[n]},y_{[n]}\}$
of $\iota^* \oplus \iota=(\mathfrak{h}^*\oplus\mathfrak{h})^\perp$. From
\begin{gather*}
\kappa^L(x_{i+1}-x_i,x_{[n]})=na_1(x_{i+1}-x_i)+nb_1(y_{i+1}-y_i)=0,
\\
\kappa^L(y_{i+1}-y_i,y_{[n]}=na_2(x_{i+1}-x_i)+nb_2(y_{i+1}-y_i)=0
\end{gather*}
we obtain $a_1=b_1=a_2=b_2=0$. We~also require that
\begin{gather*}
\kappa^L(x_{i+1}\!-x_i,y_{[n]})=(a_3\!+na_5\!-(a_5\!+a_6))(x_{i+1}\!-x_i)
\!+(b_3\!+nb_5\!-(b_5\!+b_6))(y_{i+1}\!-y_i)=0
\end{gather*}
and
\begin{gather*}
\kappa^L(x_{[n]},y_{i+1}\!-y_i)=(a_3\!+na_6\!-(a_5\!+a_6))(x_{i+1}\!-x_i)
\!+(b_3\!+nb_6\!-(b_5\!+b_6))(y_{i+1}\!-y_i)=0,
\end{gather*}
and thus all four coefficients are zero. Using the results of the image constraint
to simplify those coefficients yields
\begin{alignat*}{3}
&a_4-a_5-a_7=0,\qquad&& b_4-b_5-b_7=0,& \\
&a_4-a_6-a_7=0,\qquad && b_4-b_6-b_7=0,&
\end{alignat*}
from which it follows that $a_5=a_6$ and $b_5=b_6$.
The remaining extension requirement imposes no further constraints on the parameters because
one verifies $\kappa^L(x_{[n]},y_{[n]})=0$ using $a_5+a_6+na_7=b_5+b_6+nb_7=0$.

To analyze the PBW conditions~$\gamma_1=\gamma_2=\gamma_4=\gamma_5=0$ in~\eqref{ClearFirstObstr}
it will help to first observe that the above constraints $a_5=a_6$ and $a_5+a_6+na_7=0$
yield that
\begin{gather*}
a_5=a_6=-\frac{n}{2}a_7,
\end{gather*}
and hence that
\begin{gather*}
a_4=a_5+a_7=-\frac{n-2}{2} a_7\qquad \text{and} \qquad
a_3=-na_4=\frac{n(n-2)}{2}a_7,
\end{gather*}
with corresponding expressions in terms of $b_7$ for $b_3$, $b_4$, $b_5$, and $b_6$.
These allow the simplification
\begin{gather*}
\phi_{1,1}(x_i,x_j,y_k)
=[-b_5a_6\!+b_7(a_3\!-\!a_5\!-\!a_6)](x_i\!-\!x_j)\!+[b_5(b_5\!+nb_7)\!+b_7(b_3\!-b_5\!-b_6)](y_i\!-\!y_j)
\\ \hphantom{\phi_{1,1}(x_i,x_j,y_k)}
{}=\frac{n^2}{4}a_7b_7(x_i-x_j)+\frac{n^2}{4}b_7^2(y_i-y_j).
\end{gather*}
Similarly,
\begin{gather*}
\phi_{1,1}(y_i,y_j,x_k)=\frac{n^2}{4}a_7^2(x_i-x_j)+\frac{n^2}{4}a_7b_7(y_i-y_j).
\end{gather*}

Requiring each of these to be zero forces $a_7=b_7=0$, and thus $a_i=b_i=0$ for $3 \leq i \leq 6$
as~well. Since we already have $a_1=b_1=a_2=b_2=0$, this proves there are no
Lie orbifold algebra maps for $S_n$ acting on the doubled standard subrepresentation $\mathfrak{h}^* \oplus \mathfrak{h}$
with $\kappa^L \not\equiv 0$.
\end{proof}

For maps with linear part supported only off the identity there is instead a
three-parameter family of~Drinfeld orbifold algebra maps that generalize the
commutator relations for the rational Cherednik algebra $H_{0,c}$.

\begin{Theorem}
\label{thm:DOAMapsDoubledStd}
Let $S_n$ $(n \geq 3)$ act via the doubled standard representation on
the space $\mathfrak{h}^*\oplus\mathfrak{h} \cong \mathbb{C}^{2n-2}$ spanned by $\bar{x}_i=x_i+\frac{1}{n}x_{[n]}$ and $\bar{y}_i=y_i+\frac{1}{n}y_{[n]}$ with $1\leq i \leq n$. 
Let $a^\perp,b^\perp,c \in \mathbb{C}$.
All~Drinfeld orbifold algebra maps with
nonzero linear part supported only off the identity have the form
$\kappa^L+\kappa^C$
defined for $1 \leq i \neq j \leq n$ by
\begin{gather*}
\kappa^L(\bar{x}_i,\bar{x}_j)=\kappa^L(\bar{y}_i,\bar{y}_j)=0, \\
\kappa^L(\bar{x}_i,\bar{y}_i)=-\frac{n}{2}\sum\limits_{k\neq i} \big(a^\perp(\bar{x}_{i}+\bar{x}_{k})+b^\perp(\bar{y}_{i}+\bar{y}_{k})\big) \otimes (ik), \\
\kappa^L(\bar{x}_i,\bar{y}_j)=\frac{n}{2}\big(a^\perp (\bar{x}_{i}+\bar{x}_{j})+b^\perp(\bar{y}_{i}+\bar{y}_{j})\big) \otimes (ij)
\end{gather*}
and
\begin{alignat*}{3}
& \kappa^C(\bar{x}_i,\bar{x}_j)=\frac{n^2}{4}\big(b^\perp\big)^2 \sum\limits_{k\neq i,j} (ijk)-(kji), \quad&&
\kappa^C(\bar{y}_i,\bar{y}_j)=\frac{n^2}{4}\big(a^\perp\big)^2 \sum\limits_{k\neq i,j} (ijk)-(kji),& \\
&\kappa^C(\bar{x}_i,\bar{y}_i)=c\sum\limits_{k\neq i} (ik), \quad&&
\kappa^C(\bar{x}_i,\bar{y}_j)=-c(ij)-\frac{n^2}{4}a^\perp b^\perp\sum\limits_{k\neq i,j} (ijk)-(kji).&
\end{alignat*}
\end{Theorem}

\begin{proof}Suppose $\kappa=\kappa^L+\kappa^C$ is a Drinfeld orbifold algebra map for the doubled
standard representation with $\kappa^C\colon \bigwedge^2 (\mathfrak{h}^*\oplus\mathfrak{h}) \to \mathbb{C} S_n$ and
nonzero $\kappa^L\colon \bigwedge^2 (\mathfrak{h}^*\oplus\mathfrak{h}) \to (\mathfrak{h}^*\oplus\mathfrak{h}) \otimes \mathbb{C} S_n$ supported
only off the identity. Extend $\kappa$ as described in Corollary~\ref{prop:Extension}
to yield a Drinfeld orbifold algebra map for $S_n$ acting on $V=W^* \oplus W \cong \mathbb{C}^{2n}$
via the doubled permutation representation. By~Theorem~\ref{thm:DOAMapsDoubledPerm} the
possible forms of such extensions are
$\kappa^L_{\refl}+\kappa^C_1+\kappa^C_{\refl}+\kappa^C_{\tri}$
satisfying the PBW conditions~\eqref{Obstr2PhiC1C2C3L2=0a} and~\eqref{Obstr2PhiC1C2C3L2=0b}.

We start by imposing the condition $\im \kappa^L \subseteq \mathfrak{h}^*\oplus\mathfrak{h}$.
Use $\bar{x}_i=x_i-\tfrac{1}{n}x_{[n]}$, $\bar{y}_i=y_i-\tfrac{1}{n}y_{[n]}$,
$\bar{x}_{ij}:=\bar{x}_i+\bar{x}_j$, and $\bar{y}_{ij}:=\bar{y}_i+\bar{y}_j$ to
rewrite the nonzero values of $\kappa^L_{(ij)}$ as
\begin{align*}
\kappa^L_{(ij)}(x_i,y_i)&=\kappa^L_{(ij)}(x_j,y_j)=-\kappa^L_{(ij)}(x_i,y_j)=-\kappa^L_{(ij)}(x_j,y_i) \\
&=ax_{ij}+a^\perp x_{ij}^\perp + by_{ij}+b^\perp y_{ij}^\perp
=\big(a-a^\perp\big)x_{ij}+a^\perp x_{[n]}+\big(b-b^\perp\big)y_{ij}+b^\perp y_{[n]}
\\
&=\big(a\!-a^\perp\big)\bar{x}_{ij}\!+\frac{1}{n}\big(2a+(n\!-2)a^\perp\big) x_{[n]}\!+\big(b\!-b^\perp\big)\bar{y}_{ij}\!+\frac{1}{n}\big(2b+(n\!-2)b^\perp\big) y_{[n]}.
\end{align*}
This shows
$2a+(n-2)a^\perp=2b+(n-2)b^\perp=0$ or
\begin{equation}
\label{arelatedtoaperp}
a-a^\perp=-\frac{n}{2}a^\perp \quad \text{ and } \quad b-b^\perp=-\frac{n}{2}b^\perp
\end{equation}
and substituting these into~\eqref{Obstr2PhiC1C2C3L2=0a}
and~\eqref{Obstr2PhiC1C2C3L2=0b} yields that{\samepage
\begin{gather*}
a(\alpha-\beta)=a^\perp(\alpha-\beta)=b(\alpha-\beta)=b^\perp(\alpha-\beta)=0.
\end{gather*}
Thus since $\kappa^L \not \equiv 0$ we must also have $\alpha=\beta$ in $\kappa^C_1$.}

Next consider conditions arising from $\iota^*\oplus\iota\subseteq\ker\kappa$ in
Corollary~\ref{prop:Extension}, i.e., $\kappa^C(w,v)=\kappa^C(u,v)=0$
for all $w$ in the basis $\{x_{i+1}-x_i,\, y_{i+1}-y_i \mid 1 \leq i \leq n-1 \}$
of $\mathfrak{h}^*\oplus\mathfrak{h}$ and~$u$~$v$ in the basis $\{x_{[n]},y_{[n]}\}$
of $\iota^* \oplus \iota=(\mathfrak{h}^*\oplus\mathfrak{h})^\perp$.

For $\kappa^C_1$, since $\kappa^C_1(x_i,x_j)=\kappa^C_1(y_i,y_j)=0$ and
$\kappa^C_1(x_i,y_{[n]})=\kappa^C_1(x_{[n]},y_i)=\alpha+(n-1)\beta$ for
any $1 \leq i,j \leq n$, the only extension condition that is not
automatically satisfied is $\kappa^C_1(x_{[n]},y_{[n]})=n(\alpha+(n-1)\beta)=0$.
Together with $\alpha=\beta$ this forces $\alpha=\beta=0$ and thus $\kappa^C_1 \equiv 0$.

For $\kappa^*_{\refl}$ and $\kappa^C_{\tri}$, the definitions of $\kappa^L_g$
and $\kappa^C_g$ when $g$ is a transposition and of $\kappa^C_g$ when $g$
is a $3$-cycle yield that
\begin{equation}
\label{orbiteqn}
\kappa^*_g\bigg(v,\sum\limits_{h \in \langle g \rangle} \actby{h}v\bigg)=
\kappa^*_g\bigg(v,\sum\limits_{h \in \langle g \rangle} \actby{h}v^*\bigg)=0
\end{equation}
for all $v \in \{x_1,\dots,x_n,y_1,\dots,y_n\}$.
This in turn implies that all of the extension conditions hold
for $\kappa^*_{\refl}$ and $\kappa^C_{\tri}$, yielding no further constraints
on parameters.

We now evaluate $\kappa^L_{\refl}$, $\kappa^C_{\refl}$, and $\kappa^C_{\tri}$
at pairs of vectors
from $\{\bar{x}_1,\dots,\bar{x}_n,\bar{y}_1,\dots,\bar{y}_n\}$. For~$g$ a transposition, $i$ an index moved by $g$, $\bar{v}$
a vector in $\{\bar{x}_i,\actby{g}\bar{x}_i,\bar{y}_i,\actby{g}\bar{y}_i\}$,
and $*=L$ or $*=C$, we~have $\kappa^*_g(\bar{v},\actby{g}\bar{v})=0$
and we use~\eqref{orbiteqn} to observe that for $1 \leq i \leq n$,
\begin{gather*}
\kappa^*_g(\bar{x}_i,\bar{y}_i)=-\kappa^*_g(\bar{x}_i,\actby{g}\bar{y}_i)=
\kappa^*_g(x_i,y_i)-\frac{1}{n}\kappa^*_g(x_i,y_{[n]})
-\frac{1}{n}\kappa^*_g(x_{[n]},y_j)+\frac{1}{n^2}\kappa^*_g(x_{[n]},y_{[n]})
\\ \hphantom{\kappa^*_g(\bar{x}_i,\bar{y}_i)}
{}=\kappa^*_g(x_i,y_i).
\end{gather*}
It then follows from~\eqref{arelatedtoaperp} that
\begin{gather*}
\kappa^L_g(\bar{x}_i,\bar{y}_i)=-\kappa^L_g(\bar{x}_i,\actby{g}\bar{y}_i)
=-\frac{n}{2}a^\perp \bar{x}_{ij}-\frac{n}{2}b^\perp \bar{y}_{ij},
\\
\kappa^C_g(\bar{x}_i,\bar{y}_i)=-\kappa^C_g(\bar{x}_i,\actby{g}\bar{y}_i)=c.
\end{gather*}
For $g$ a $3$-cycle, by the orbit property in~\eqref{orbiteqn} and by~\eqref{arelatedtoaperp}
we see that
\begin{gather*}
\kappa^C_g(\bar{v},\bar{v}^*)=0, \\
\kappa^C_g(\bar{v},\actby{g}\bar{v})=
\begin{cases}
 	\dfrac{n^2}{4}\big(a^\perp\big)^2 & \text{if}\quad v \in W, \\[2ex]
 	\dfrac{n^2}{4}\big(b^\perp\big)^2 & \text{if}\quad v \in W^*,
\end{cases}
\\
\kappa^C_g(\actbylessspace{g}\bar{v},\bar{v}^*)=-\kappa^C_g(\bar{v},\actbylessspace{g}\bar{v}^*)=
\frac{n^2}{4}a^\perp b^\perp.
\end{gather*}
These components produce the given definition of $\kappa^L+\kappa^C$.
\end{proof}

\begin{Remark}\label{RecoverGHAs}
In the case $\kappa^L_1=\kappa^L_{\refl} \equiv 0$ then also $\kappa^C_{\tri}\equiv 0$,
$\alpha=-(n-1)\beta$, and the restriction of the constant $2$-cochain
$\kappa^C=\kappa^C_1+\kappa^C_{\refl}$ to $\bigwedge^2 (\mathfrak{h}^*\oplus\mathfrak{h})$ is given by
\begin{alignat*}{3}
\kappa^C_1(\bar{x}_i,\bar{x}_j)&=\kappa^C_1(\bar{y}_i,\bar{y}_j)=0, &\qquad \kappa^C_1(\bar{x}_i,\bar{y}_i)&=-(n-1)\beta, &\qquad \kappa^C_1(\bar{x}_i,\bar{y}_j)&=\beta,
\\
\kappa^C_g(\bar{x}_i,\bar{x}_j)&=\kappa^C_g(\bar{y}_i,\bar{y}_j)=0, &\qquad \kappa^C_g(\bar{x}_i,\bar{y}_i)&=c, &\qquad \kappa^C_g(\bar{x}_i,\bar{y}_j)&=-c,
\end{alignat*}
where $\beta,c \in \mathbb{C}$, $g$ is a transposition, and $1 \leq i \neq j \leq n$.
This corresponds to the rational Cherednik algebra $H_{n\beta,c}$.
\end{Remark}

In the theory of rational Cherednik algebras, $H_{t,c}$, for the
symmetric group, a natural isomorphism between $H_{t,c}$ and
$H_{\lambda t, \lambda c}$ when $\lambda \in \mathbb{C}^\times$
means that only two distinct cases need be considered, $t\neq 0$ and $t=0$.
Theorems~\ref{thm:LOAMapsDoubledStd} and~\ref{thm:DOAMapsDoubledStd}
show that in the first case there are no further deformations in polynomial
degree one with the linear part of the parameter map supported only on the
identity while there is a three-parameter family of such deformations in the
second case with the linear part of the parameter map supported only off the
identity.

\begin{Remark}\label{NoExtensionsOfExamples}
What if $\kappa^L$ were supported both on and off the
identity? The specializations of~parameter values for $\kappa^L_1$ in
part (1) and for $\kappa^L_{\refl}$ in part (3) of
Theorem~\ref{thm:DOAMapsCombined} on the combined lift of
$\kappa^L_1+\kappa^L_{\refl}$ means that
\begin{gather*}
\kappa^L_1(x_i,y_i)=\kappa^L_1(x_i,y_j)=a_4x_{[n]}+b_4y_{[n]},
\\
\kappa^L_{(ij)}(x_i,y_i)=-\kappa^L_{(ij)}(x_i,y_j)=
-\frac{n}{2}a^\perp(\bar{x}_i+\bar{x}_j)
-\frac{n}{2}b^\perp(\bar{y}_i+\bar{y}_j).
\end{gather*}
But then $\im\big(\kappa^L_1+\kappa^L_{\refl}\big) \subseteq \mathfrak{h}^*\oplus\mathfrak{h}$ would
require $a_4=b_4=0$ so $\kappa^L_1\equiv 0$. This combined with part (2) of
Theorem~\ref{thm:DOAMapsCombined} shows there is no Drinfeld orbifold
algebra map for $S_n$ on $\mathfrak{h}^*\oplus\mathfrak{h}$ with linear part supported both on
and off the identity which extends to a map $\kappa$ of the form in
Theorem~\ref{thm:DOAMapsCombined}. But since Theorem~\ref{thm:DOAMapsCombined}
is not exhaustive, it is not clear whether there exist such maps in general.
\end{Remark}

\section{Descriptions of degree-one rational Cherednik algebras}\label{Deg1RCAs}


Here we present, via generators and relations, degree-one PBW deformations
of the skew group algebras $S(W^* \oplus W) \# S_n$ and $S(\mathfrak{h}^* \oplus \mathfrak{h})\# S_n$
that result from Theorems~\ref{thm:LOAMapsDoubledPerm},~\ref{thm:DOAMapsDoubledPerm},
and~\ref{thm:DOAMapsDoubledStd} when \mbox{$n \geq 3$}. This facilitates comparison with degree-zero
deformations (i.e., rational Cherednik algebras) and with the PBW deformations of
$S(W)\# S_n$ in~\cite{FGK}. The classifications are summarized in
Tables~\ref{SummaryDOAMapsDoubledPerm} and~\ref{SummaryDOAMapsDoubledStd}. We~reiterate that the case when $n=2$ can be analyzed in similar fashion, but involves
some differences in the dimensions of spaces of pre-Drinfeld
orbifold algebra maps and in the parameter relations required in order to satisfy the
PBW conditions.

\subsection{Algebras for the doubled permutation representation}\label{AlgsDoubledPermRep}

First, the Lie orbifold algebra maps involving $17$ parameters classified in
Theorem~\ref{thm:LOAMapsDoubledPerm} yield a~variety controlling
the Lie orbifold algebras that deform $S(W^* \oplus W) \# S_n$ in degree one.
Based on representative calculations in Macaulay2~\cite{M2} we conjecture
that this projective variety is of~dimension seven. Some subvarieties of potential
interest are indicated in Table~\ref{table:subvars} in~Section~\ref{InvLieBrackets}.
When $\kappa^L_1 \equiv 0$ these Lie orbifold algebras
specialize to rational Cherednik algebras corresponding to the parameter $c$ and the general
$G$-invariant skew-symmetric bilinear form $\kappa^C_1$ involving $\alpha$ and $\beta$
(because $W$ is decomposable~--- see~\cite[proof of Theorem 1.3]{EtingofGinzburg2002}).

\begin{Theorem}[{Lie orbifold algebras for doubled permutation representation over $\mathbb{C}[t]$}]
\label{Lieorbifold}
Let $S_n$ $(n \geq 3)$ act on $V=W^* \oplus W$ with basis
$\mathcal{B}=\{x_1,\dots,x_n,y_1,\dots,y_n\}$ by the doubled permutation
representation. For~$a_1, \dots, a_7, b_1, \dots, b_7, \alpha, \beta, c \in \mathbb{C}$
subject to conditions~\eqref{ClearFirstObstr},~\eqref{ClearSecondObstr1C}, and~\eqref{ClearSecondObstr1Ref}, define
$\kappa^L=\kappa^L_1$ and $\kappa^C=\kappa^C_1+\kappa^C_{\refl}$ to be the linear and
constant cochains such that for $1\leq i\neq j \leq n$,
\begin{alignat*}{3}
&\kappa^L(x_i,x_j)=(a_1(x_i-x_j)+b_1(y_i-y_j)), \qquad&& \kappa^C(x_i,x_j)=0,& \\
& \kappa^L(y_i,y_j)=(a_2(x_i-x_j)+b_2(y_i-y_j)), \qquad &&\kappa^C(y_i,y_j)=0, &\\
&\kappa^L(x_i,y_i)=(a_3x_i+a_4\n{x}+b_3y_i+b_4\n{y}), \qquad&&
\kappa^C(x_i,y_i)=\alpha+c\sum\limits_{k \neq i} (ik),& \\
&\kappa^L(x_i,y_j)=(a_5x_i+a_6x_j+a_7\n{x}+b_5y_i+b_6y_j+b_7\n{y}) ,\qquad&&
\kappa^C(x_i,y_j)=\beta-c(ij).&
\end{alignat*}
Then the quotient $\mathcal{H}_{\kappa,t}$ of $T(V) \# S_n[t]$ by the ideal generated
by
\begin{gather*}
\big\{uv-vu - \kappa^L(u,v) t -\kappa^C(u,v) t^2 \mid u,v \in \mathcal{B} \big\}
\end{gather*}
is a Lie orbifold algebra over $\mathbb{C}[t]$.
In fact, the algebras $\mathcal{H}_{\kappa,1}$ are precisely the Drinfeld orbifold
algebras such that $\kappa^L$ is supported only on the identity.
\end{Theorem}

\begin{table}[h!]
\begin{center}
\caption{Classification of~Drinfeld orbifold algebra maps for $S_n$ acting on $W^*\oplus W$.}\label{SummaryDOAMapsDoubledPerm}

\vspace{2mm}
\renewcommand{\arraystretch}{1.2}%
\begin{tabular}{l|l|c|l}
\hline
\multicolumn{1}{c|}{Linear part $\kappa^L$} & \multicolumn{1}{c|}{Constant part $\kappa^C$} & \multicolumn{1}{c|}{Parameter relations} & \multicolumn{1}{c}{Reference}
\\
\hline
$\kappa_1^L$ & 0 & \eqref{ClearFirstObstr} & Theorem~\ref{thm:LOAMapsDoubledPerm}
\\
& $\kappa_1^C$ & \eqref{ClearFirstObstr}--\eqref{ClearSecondObstr1C}
\\
& $\kappa_1^C +\kappa_{\refl}^C$ with $\kappa_{\refl}^C\not\equiv0$ & \eqref{Obstr2kappaCrefc}--\eqref{Obstr2kappaC1bSimplified}
\\
\hline
$\kappa_{\refl}^L$ & $\kappa_{\tri}^C$ & none & Theorem~\ref{thm:DOAMapsDoubledPerm}
\\
& $\kappa_{\tri}^C+\kappa_{\refl}^C$ & none
\\
& $\kappa_{\tri}^C+\kappa_{\refl}^C+\kappa_1^C$ &\eqref{Obstr2PhiC1C2C3L2=0a}--\eqref{Obstr2PhiC1C2C3L2=0b}
\\
\hline
$\kappa_1^L+\kappa_{\refl}^L$ & $\kappa_{\tri}^C+\kappa_{\refl}^C+\kappa_1^C$ & Theorem~\ref{thm:DOAMapsCombined}(1)--(3) & Theorem~\ref{thm:DOAMapsCombined}
\\
& ? & ?
\\
\hline
0 & $\kappa_1^C+\kappa_{\refl}^C$ & \text{none}
\\
\hline
\end{tabular}
\end{center}\vspace{-1ex}
{\small When the indicated parameter relations are satisfied, the map $\kappa=\kappa^L+\kappa^C$ is a Drinfeld orbifold algebra map.
The question marks indicate there could be further maps with $\kappa^L=\kappa_1^L+\kappa_{\refl}^L$.}
\end{table}

Second, for $\kappa^L$ supported only off the identity, Theorem~\ref{thm:DOAMapsDoubledPerm}
shows that by comparison there is only a seven-parameter family of~Drinfeld orbifold algebra
maps and these are controlled by a projective variety which, according to a few representative
calculations in Macaulay2~\cite{M2}, appears to be four-dimensional. The resulting algebras
also specialize to rational Cherednik algebras
parametrized by $\alpha$, $\beta$, and $c$ when $\kappa^L_{\refl}=\kappa^C_{\tri} \equiv 0$.

\begin{Theorem}[{Drinfeld orbifold algebras for doubled permutation representation over $\mathbb{C}[t]$}]
\label{permexamples}
Let $S_n$ $(n \geq 3)$ act on $V=W^* \oplus W$ with basis
$\mathcal{B}=\{x_1,\dots,x_n,y_1,\dots,y_n\}$ by the doubled permutation
representation. Suppose $a,a^\perp,b,b^\perp,c,\alpha,\beta \in \mathbb{C}$ 
satisfy conditions~\eqref{Obstr2PhiC1C2C3L2=0a} and~\eqref{Obstr2PhiC1C2C3L2=0b}.
Define $\kappa^L=\kappa_1^L$ and $\kappa^C=\kappa_1^C+\kappa_{\refl}^C+\kappa_{\tri}^C$
to be the cochains such that for $1 \leq i \neq j \leq n$,
\begin{gather*}
\kappa^L(x_i,x_j)=\kappa^L(y_i,y_j)=0, \\
\kappa^L(x_i,y_i)=\textstyle\sum\limits_{k \neq i}\big( \big(a-a^\perp\big)x_{i,k}+a^\perp x_{[n]}+\big(b-b^\perp\big)y_{i,k}+b^\perp y_{[n]}\big) \otimes (ik), \\
\kappa^L(x_i,y_j)=-\big(\big(a-a^\perp\big)x_{i,j}+a^\perp x_{[n]}+\big(b-b^\perp\big)y_{i,j}+b^\perp y_{[n]}\big) \otimes (ij)
\end{gather*}
and
\begin{gather*}
\kappa^C(x_i,y_j)=\beta-c(ij)-\big(a-a^\perp\big)\big(b-b^\perp\big)\textstyle\sum\limits_{k \neq i,j} (ijk)-(kji), \\
\kappa^C(x_i,x_j)=\big(b-b^\perp\big)^2 \textstyle\sum\limits_{k \neq i,j} (ijk)-(kji), \\
\kappa^C(y_i,y_j)=\big(a-a^\perp\big)^2 \textstyle\sum\limits_{k \neq i,j} (ijk)-(kji), \\
\kappa^C(x_i,y_i)=\alpha+c\textstyle\sum\limits_{k \neq i} (ik).
\end{gather*}
Then the quotient $\mathcal{H}_{\kappa,t}$ of $T(V) \# S_n[t]$ by the ideal generated by
\begin{gather*}
\big\{uv-vu - \kappa^L(u,v) t -\kappa^C(u,v) t^2 \mid u,v \in \mathcal{B} \big\}
\end{gather*}
is a Drinfeld orbifold algebra over $\mathbb{C}[t]$.
Further, the algebras $\mathcal{H}_{\kappa,1}$ are precisely the Drinfeld orbifold
algebras such that $\im\kappa^L_g\subseteq V^g$ for each $g\in S_n$ and $\kappa^L$
is supported only off the identity.
\end{Theorem}

An analogous statement may be made for algebras constructed from the family of
lifts of $\kappa^L_1+\kappa^L_{\refl}$ described in Theorem~\ref{thm:DOAMapsCombined}
but is omitted here.

\begin{table}[h!]
\begin{center}
\caption{Classification of~Drinfeld orbifold algebra maps for $S_n$ acting on $\mathfrak{h}^* \oplus \mathfrak{h}$.}\label{SummaryDOAMapsDoubledStd}

\vspace{2mm}

\renewcommand{\arraystretch}{1.2}%
\begin{tabular}{l|l|c|l}
\hline
\multicolumn{1}{c|}{Linear part $\kappa^L$} & \multicolumn{1}{c|}{Constant part $\kappa^C$} & \multicolumn{1}{c|}{Parameter relations} & \multicolumn{1}{c}{Reference}
\\
\hline
$\kappa_{\refl}^L$ & $\kappa_{\tri}^C$ & $2a+(n-2)a^{\perp}=0$ & Theorem~\ref{thm:DOAMapsDoubledStd}
\\
 & & $2b+(n-2)b^{\perp}=0$ &
 \\
& $\kappa_{\tri}^C+\kappa_{\refl}^C$ & $2a+(n-2)a^{\perp}=0$ &
\\
& & $2b+(n-2)b^{\perp}=0$ &
\\
\hline
0 & $\kappa_{\refl}^C$ & \text{none} & Remark~\ref{RecoverGHAs} \\
& $\kappa_1^C$ & $\alpha+(n-1)\beta=0$ & \\
& $\kappa_1^C+\kappa_{\refl}^C$ & $\alpha+(n-1)\beta=0$ &
\\
\hline
\end{tabular}
\end{center}\vspace{-1ex}
{\small When the parameter relations hold, the map $\kappa=\kappa^L+\kappa^C$ is a Drinfeld orbifold algebra map.}
\end{table}

\subsection{Algebras for the doubled standard representation}\label{AlgsDoubledStdRep}

By Theorem~\ref{thm:LOAMapsDoubledStd}
the only Lie orbifold algebras for $S_n$ acting on $\mathfrak{h}^* \oplus \mathfrak{h}$ by the doubled
standard representation are the known rational Cherednik algebras $H_{n\beta,c}$.
However, by Theorem~\ref{thm:DOAMapsDoubledStd}
there is in this case a three-parameter family of~Drinfeld orbifold algebras which are
not graded Hecke algebras, but which specialize when $a^\perp=b^\perp=0$ to the rational
Cherednik algebras $H_{0,c}$ for $S_n$.

\begin{Theorem}[{Drinfeld orbifold algebras for doubled standard representation over $\mathbb{C}[t]$}]\label{refexamples}
Let~$S_n$ $(n \geq 3)$ act on $\mathfrak{h}^* \oplus \mathfrak{h}$
by the doubled standard representation. For~$a^\perp,b^\perp,c \in \mathbb{C}$ and
$\bar{x}_i$ and $\bar{y}_i$ as in~\eqref{xbar} define
$\kappa^L$ and $\kappa^C$ as in Theorem~$\ref{thm:DOAMapsDoubledStd}$.
Then the quotient $\mathcal{H}_{\kappa,t}$ of $T(\mathfrak{h}^* \oplus \mathfrak{h}) \# S_n[t]$ by the ideal generated
by
\begin{gather*}
\big\{\bar{u}\bar{v}-\bar{v}\bar{u}-\kappa^L(\bar{u},\bar{v})t-\kappa^C(\bar{u},\bar{v})t^2 \mid \bar{u}, \bar{v} \in \bar{\mathcal{B}}\big\}
\end{gather*}
is a Drinfeld orbifold algebra of $S(\mathfrak{h}^* \oplus \mathfrak{h})\# S_n$ over $\mathbb{C}[t]$.
Further, the algebras $\mathcal{H}_{\kappa,1}$ are precisely the Drinfeld orbifold
algebras such that $\im\kappa^L_g\subseteq (\mathfrak{h}^*\oplus \mathfrak{h})^g$ for each $g\in S_n$ and $\kappa^L$
is supported only off the identity.
Specializing $a^\perp=b^\perp=0$ yields the rational Cherednik algebra~$H_{0,c}$.
\end{Theorem}

\subsection*{Acknowledgements}
We thank the referees for helpful suggestions and questions that improved the writing of the paper, particularly those that prompted us to add Proposition~\ref{prop:directsum}, improve Section~\ref{algvars}, and reorganize overall. We~also thank Lily Silverstein for helpful conversations related to Section~\ref{algvars}.

\pdfbookmark[1]{References}{ref}
\LastPageEnding

\end{document}